\documentclass[a4paper, leqno, 11pt]{amsart}
\usepackage{amsfonts}
\usepackage{amssymb}
\usepackage{graphicx}
\graphicspath{ {CePC/Bureau/} }
\usepackage[applemac,utf8]{inputenc}
\usepackage{color}
\usepackage{hyperref}\hypersetup{pdfborder=0 0 0}
\usepackage{mathdots}
\usepackage{float}
\newtheorem{lemma}{Lemma}[section]

\newtheorem{definition}[lemma]{Definition}
\newtheorem{proposition}[lemma]{Proposition}
\newtheorem{theorem}[lemma]{Theorem}
\newtheorem{remark}[lemma]{Remark}

\newtheorem{counterexample}[lemma]{Counterexample}
\newtheorem{corollary}[lemma]{Corollary}

\title[Infinite products of nonnegative matrices]{Normalized image of a vector by an infinite product of nonnegative matrices}
\author[A. Thomas]{Alain Thomas}
\address{Alain Thomas,
448 all\'ee des Cantons, 83640 Plan d'Aups Sainte Baume, France}
\email{alain-yves.thomas@laposte.net}
\keywords{Infinite products of nonnegative matrices; Projective spaces; Sofic measures; Weak Gibbs measures; %Functions of finite Markov chains; 
Bernoulli convolutions}
\subjclass{15B48, 28A12}

\date{June 2024}

\begin{document}

\maketitle

\centerline{\it Dedicated to the memory of Adrien Douady}

\begin{abstract}To prove that a measure, linearly representable by means of a finite set of nonnegative matrices $\mathcal M$, has the weak-Gibbs property, one check the uniform convergence (on $\mathcal M^\mathbb N$) of the sequence of vectors $\frac{A_1\cdots A_nc}{\Vert A_1\cdots A_nc\Vert}$ ($c$ positive column-vector). The main theorem gives a sufficient condition for this sequence to converge pointwise. This theorem generalizes the Birkhoff contraction method because it can be used even if the matrices have many zero entries. We also look at the convergence of the sequence of matrices $\frac{A_1\cdots A_n}{\Vert A_1\cdots A_n\Vert}$. The measures defined by Bernoulli convolution are in certain cases linearly representable; we give two example of weak-Gibbs Bernoullt convolutions, by using the Birkhoff contraction coefficient for the first and the theorem for the second. Furthermore we explicit the relationship between the notions of Bernoulli convolution, fundamental curves and lattice two-scale difference equations.\end{abstract}

\section{Introduction}

By Carathéodory's extension theorem one can define as follows a probability measure~$\nu$ on the Borel algebra of the set $\Omega_a:=\{0,1,\dots,a-1\}^\mathbb N$. Given some nonnegative row-vectors $r_0,\dots,r_{a-1}$ (i.e. all their entries are nonegative), some square nonnegative matrices $M_0,\dots,M_{a-1}$ and a positive column-vector $c$ (i.e. all the entries of $c$ are positive) such that $\sum_ir_ic=1$ and $\sum_iM_ic=c$, 
 we set\begin{equation}\label{lin}\begin{array}{l}\nu([\omega_1\dots\omega_n]):=r_{\omega_1}M_{\omega_2}\cdots M_{\omega_n}c\\\text{for any cylinder-set }[\omega_1\dots\omega_n]:=\{\xi\in\Omega_a\;;\;\xi_1=\omega_1,\dots,\xi_n=\omega_n\}.\end{array}\end{equation}

The measure defined in (\ref{lin}) is called "linearly representable" \cite[Definition~4.12]{BP}, or "sofic", "hidden Markov chain", "function of Markov chain", "probabilistic function of Markov chain" (\cite[\S 2.4 and Theorem~4.20]{BP}, \cite{Eri,Ley,Petrie}, \cite[\S2]{T2}).

\begin{remark}The definition (\ref{lin}) generalizes the one of a Markov probability measure $\mu$ by the classical formula $\mu[\omega_1\dots\omega_n]:=p_{\omega_1}p_{\omega_1,\omega_2}\cdots p_{\omega_{n-1},\omega_n}$: this is the special case $r_k:=\left(\begin{smallmatrix}0&\dots&p_k&\dots&0\end{smallmatrix}\right)$, $c=\left(\begin{smallmatrix}1\\\vdots\\1\end{smallmatrix}\right)$, $M_k:=\left(\begin{smallmatrix}0&\dots&p_{0,k}&\dots&0\\\vdots&\ddots&\vdots&\ddots&\vdots\\0&\dots&p_{a-1,k}&\dots&0\end{smallmatrix}\right)$, where the  columns of $r_k$ and $M_k$ are null except the $(k+1)^{\rm th}$).\end{remark}

\begin{remark} To prove that a linearly representable measure has the weak-Gibbs property (see \cite[Definition 1.1]{FO}, \cite[Definition 1.1]{IY}, \cite[\S5]{Y}), one checks for any positive vector $c$ that the sequence$$\textstyle c_{\omega,n}:=\frac{M_{\omega_1}\cdots M_{\omega_n}c}{\Vert M_{\omega_1}\cdots M_{\omega_n}c\Vert}$$converges uniformly in $(\omega_n)_{n\in\mathbb N}\in\{0,\dots,a-1\}^\mathbb N$ (see \cite[Lemma~1.2]{FO} and Theorem \ref{weak-Gibbs} below). 

The Markov probability measures are a special case of weak-Gibbs measures: the potential involved in Definition \ref{2definitions}(i) is $\phi(\omega)=\log p_{\omega_1,\omega_2}$.\end{remark}

For any sequence of nonnegative $d\times d$ matrices $(A_n)_{n\in\mathbb N}$ we use the notation$$\textstyle P_n:=A_1\cdots A_n\text{ and }P_{m,n}:=A_{m+1}\cdots A_n\ (n>m\ge0),P_0=P_{n,n}:=I_d.$$It seems that in many cases the sequence of vectors $\frac{P_nc}{\Vert P_nc\Vert}$ converges while the sequence of matrices $\frac{P_n}{\Vert P_n\Vert}$ diverges (\S1 and \S4 respectively). The convergence of $\frac{P_nc}{\Vert P_nc\Vert}$ can be proved by different methods. We cannot use the known results about the infinite products of matrices \cite[\dots]{Koz13,BL85,F,M,M+} for the following reason:  the sequence $P_n$ either diverges or converges to $0$ if the set $\mathcal M_\infty:=\{M\;;\;A_n=M\text{ for infinitely many }n\}$ is finite and if the elements of this set do not have a common left-eigenvector for the eigenvalue 1. This result is a direct consequence of \cite[Theorem~1]{EF97}. So, if the matrices $A_n$ belong to a finite set, we are restricted to suppose that the matrices of $\mathcal M_\infty$ have a common left-eigenvector for the eigenvalue 1, to write $\lim_{n\to\infty}\frac{P_nc}{\Vert P_nc\Vert}=\frac{Ac}{\Vert Ac\Vert}$ if $A=\lim_{n\to\infty}P_n$ exists.

The Birkhoff contraction coefficient computed in \cite[Theorems 3.12 and 3.10]{Sen81}, is less than 1 when the matrices $A_n$ are positive, implying the convergence of the sequence $\frac{P_nc}{\Vert P_nc\Vert}$. If the matrices $A_n$ have many null entries, it is appropriate to seek another explanation for the convergence of $\frac{P_nc}{\Vert P_nc\Vert}$. The method we use in Theorem~\ref{rankone} and its application \S\ref{appli} is to write $P_n$ as a product of product-matrices:$$\textstyle P_n=P_{s_1}P_{s_1,s_2}\cdots P_{s_{k-1},s_k}P_{s_k,n}\ (s_{k+1}\le n<s_{k+2}),$$choosing the differences $s_{i+1}-s_i$ large enough  for the product-matrices $P_{s_i,s_{i+1}}$ ($1\le i\le k-1$) and $P_{s_k,n}$ to have some suitable properties (hypotheses $(H_1)$ and $(H_2)$ of the theorem) that imply the convergence of $\frac{P_nc}{\Vert P_nc\Vert}$. %The proof of the following theorem is different from the proof by the Birkhoff contraction property but leads to the same conclusion, the uniqueness of the limit-vector.

Section \ref{functional} concerns some equations whose unknown is a function (two-scale difference equations) or a probability measure (Bernoulli convolutions). One uses the infinite products of matrices to solve the equation in the first case, or to study the properties of the probability measure (which always exist) in the second case. The two-scale difference equations have the form $f(x)=\sum_{n=0}^Nc_nf(kx-n)$ with $1<k\in\mathbb R$ \cite{DL91,DL92,DL92Two,BeP} and, when such an equation has a solution $f\in L^1$ with nonnull integral, one has (by integration) $\sum_{n=0}^Np_n=1$ with $p_n=\frac{c_n}k$. If $f$ and the $c_n$ are nonnegative, the measure $\mu$ of density $f$ is a solution of the equation $\mu(B)=\sum_{n=0}^Np_n\mu(kB-n)$ for any $B$ Borel set (we are rather interested to the cases where the solution of this equation does not have a density). Given some nonnegative reals $p_0,\dots,p_N$ of sum 1, this equation always has a unique solution among the finite Borel probability measures, according to the introduction of~\cite{H}. We call it the "Bernoulli convolution $\nu_{k,p}$ associated to $k$ and to the probability vector $p=(p_0,\dots,p_N)$", by extension of the definition in \cite{Erd39,PSS}. 

The measure $\nu_{k,p}$ is linearly representable if $k\in\mathbb N$ (first examples of \S\ref{examples}), as well as in the second example, where $k\not\in\mathbb N$ is chosen among the Pisot numbers. We specify their weak-Gibbs properties, by using the Birkhoff contraction method in the first example and Theorem \ref{rankone} in the second.

Section \ref{limpoints}, about the divergence of the sequence of matrices $\frac{P_n}{\Vert P_n\Vert}$, is independent of the previous sections. We first deduce from \cite[Theorem~1.1]{F} (by applying this theorem to the transpose of the matrix $B_n=\frac{A_n}{\Vert A_n\Vert}$) that, if the matrices $A_n$ are positive and if the sequence $\frac{A_n}{\Vert A_n\Vert}$ converges to a positive matrix, then the sequence $\frac{P_n}{\Vert P_n\Vert}$ converges to a rank-one matrix. We give a partial converse of this result in Theorem~\ref{aediv}: for any sequence $(A_n)_{n\in\mathbb N}$ of complex-valued matrices such that the sequence $\frac{P_n}{\Vert P_n\Vert}$ converges, either the sequence $\frac{A_n}{\Vert A_n\Vert}$ converges or this sequence has several limit-points with a common left-eigenvector.

In Appendix \ref{mf} we give the definition of the weak-Gibbs property and we explain how it can be used to specify the local properties of a probability measure; we give the proof (Theorem \ref{weak-Gibbs}) that, for any linearly representable measure, the uniform convergence of the sequence $c_{\omega,n}$ implies the weak-Gibbs property. In Appendices \ref{proofcvandGibbs} and \ref{proofbetareal} we give the proofs relating to the examples of Section \ref{examples}.

Some general results about the convergence of the infinite products of matrices can be found in the introductions of \cite{DL92Two} and \cite{BW}:

* the RCP sets are the sets of matrices $\{M_k\}_{k\in K}$ such that the limit-matrix $M(\omega):=\lim_{n\to\infty}M_{\omega_1}\cdots M_{\omega_n}$ exists for any sequence $\omega\in K^\mathbb N$,

* the continuous RCP sets are the finite sets $\{M_0,\dots,M_{b-1}\}$ such that $M(\cdot)$ is continuous for the usual topology on $\Omega_b$,

* the real-continuous RCP sets are the finite sets $\{M_0,\dots,M_{b-1}\}$ for which $M(\omega(\cdot))$ is continuous, where $\omega(x)=(\omega_n)_{n\in\mathbb N}$ denotes the expansion of $x\in[0,1)$ in base $b$ ($x=\sum_{n\in\mathbb N}\omega_n/b^n$).

The sixth section of \cite{DL92Two} give some applications of the RCP~sets: nonhomogeneus Markov chains, self-similar curves \cite{B} like the ones of von Koch and de Rham, lattice two-scale difference equations \cite{DL91,DL92,A,X}. Other results can be found in \cite{BL85,Bo,DL01,EF97,F,Koz13,L89}. See \cite{BO,BP,Dh,Eri,Ley,Petrie} for the sofic measures and \cite{Erd39,PSS,V} for the Bernoulli convolutions. For the multifractal analysis and Lyapunov exponents, see \cite{FeandLau02,Fen03,Fen03bis,Fen04,Fen09,F24,FH10,FO,Led84,DL94,Go}.

About the special case of the $2\times2$ nonnegative matrices, the C.N.S. for the sequence $\frac{P_nc}{\Vert P_nc\Vert}$ to converge (resp. to converge uniformly), when $c$ is a positive vector and $\{M\;;\;\exists n\ A_n=M\}$ is finite, is given in \cite{T} (resp. \cite{OT13}). On the other hand, Mukherjea, Nakassis and Ratti \cite{M,M+} study the weak limit $\lambda$ of the distribution of $Y_n\cdots Y_1$ when the random matrices $Y_n$ are i.i.d in a finite set of stochastic matrices $\{M_1,\dots,M_N\}$ with $\forall i\ M_i=\left(\begin{smallmatrix}x_i&1-x_i\\y_i&1-y_i\end{smallmatrix}\right)$ and $x_i,y_i\in(0,1)$. In other words they study the random variable$$\textstyle\{1,\dots,N\}^\mathbb N\ni(k_n)_{n\in\mathbb N}\mapsto y_{k_1}+\sum_{n=2}^\infty y_{k_n}\prod_{i=1}^{n-1}\det M_{k_i}.$$In the special case where $x_{i+1}-x_i$ and $x_i-y_i$ do not depend on $i$, $
\lambda$ is a Bernoulli convolution, up to an affine transformation. From \cite[Proposition 1]{M} the density of the Bernoulli convolution $\nu_{k,p}$ associated with $k=\root m\of N$ and the probability vector $p=(\frac1N,\dots,\frac1N)$ is a piecewise polynomial of degree not exceeding $m$.

{\bf Acknowledgement. --} {\it We thank several mathematicians for their comments and help. Adrien Douady was interested to study some examples with specific methods. Eric Olivier contributes to the applications of the infinite products of matrices to the multifractal analysis and to the study of the Bernoulli convolutions. Various referees and editors give us many information about related subjects.}

\section{Rank one property, and the sequence $\frac{P_nc}{\Vert P_nc\Vert}$}\label{tsf}We use the following notations:

$u_1,\dots,u_d$ are the canonical $d$-dimensional column-vectors, and $u:=\sum_iu_i$;

$\Vert X\Vert$ is the sum of the absolute values of the entries, in the matrix or vector~$X$; in this paper, the choice of the norm does not matter because, given two norms $\mathcal N_1$ and $\mathcal N_2$, if a sequence of the form $\frac{x_n}{\mathcal N_1(x_n)}$ converges to a limit $\ell$, then the sequence $\frac{x_n}{\mathcal N_2(x_n)}$ converges to $\frac\ell{\mathcal N_2(\ell)}$;

the shift $\sigma$ is defined for any sequence $s=(s_n)_{n\in\mathbb N}$ by $\sigma(s):=(s_{n+1})_{n\in\mathbb N}$.

\subsection{Sufficient condition for the sequence $\frac{P_nc}{\Vert P_nc\Vert}$ to converge}\label{suf}

%We denote by $Q_k$ (resp. $R_n$) the matrix $P_{s_k,s_{k+1}}$ (resp. the matrix $P_{s_k,n}$, where $k$ is the integer  such that $s_{k+1}\le n<s_{k+2}$).
One can always chose an increasing sequence $(s_k)_{\in\mathbb N}$ for $\frac{P_{s_k}}{\Vert P_{s_k}\Vert}$ to converge but, when one makes the product $P_{s_k}P_{s_k,n}=P_n$, the dominant columns (in norm) of $P_{s_k}$ may disappear if the rows of $P_{s_k,n}$ with same indices as them are null, or if they are small (when $k\to\infty$) with respect to $\Vert P_{s_k,n}\Vert$. So we make two hypotheses in the following theorem.

\begin{theorem}\label{rankone}Let $\mathcal A=(A)_{n\in\mathbb N}$ be a sequence of nonnegative $d\times d$ matrices such that $\forall n\ P_n\ne0$, and suppose there exists a limit-point $P^{(s)}:=\lim_{k\to\infty}\frac{P_{s_k}}{\Vert P_{s_k}\Vert}$ such that

$(H_1)$ $\textstyle\inf\big\{\frac{\Vert P_n\Vert}{\Vert P_{s_k}\Vert\Vert P_{s_k,n}\Vert}\;;\;k,n\text{ such that }s_{k+1}\le n<s_{k+2}\big\}>0$;

$(H_2)$ setting $J=\{j\;;\;P^{(s)}u_j\ne0\}$, there exists a limit-point $Q=(q_{i,j})_{i,j}$ of the sequence $k\mapsto\frac{P_{s_k,s_{k+1}}}{\Vert P_{s_k,s_{k+1}}\Vert}$ such that$$\textstyle\forall j,j'\in J\ \exists i\in J\ q_{i,j}q_{i,j'}\ne0.$$

Then the sequence $\frac{P_nc}{\Vert P_nc\Vert}$ converges for any positive column-vector $c$, and more generally for any  $c\in\mathcal V_\mathcal A:=\{v\text{ nonnegative}\;;\;\inf_n\frac{\Vert P_nv\Vert}{\Vert P_n\Vert}>0\}$. The limit does not depend on $c$.\end{theorem}

\begin{proof}[\bf Proof]Let $\mathcal A=(A)_{n\in\mathbb N}$ and $P^{(s)}$ satisfy the hypotheses of Theorem \ref{rankone}, and let $\varepsilon_k:=\big\Vert\frac{P_{s_k}}{\Vert P_{s_k}\Vert}-P^{(s)}\big\Vert$. Using the classical inequality $\big\Vert\frac x{\Vert x\Vert}-\frac y{\Vert y\Vert}\big\Vert\le\frac{2\Vert x-y\Vert}{\Vert x\Vert}$,\begin{equation}\label{PnPs}\begin{array}{rcl}\Big\Vert\frac{P_nc}{\Vert P_nc\Vert}-\frac{P^{(s)}P_{s_k,n}c}{\Vert P^{(s)}P_{s_k,n}c\Vert}\Big\Vert&=&\Big\Vert\frac{P_nc/\Vert P_{s_k}\Vert}{\Vert P_nc\Vert/\Vert P_{s_k}\Vert}-\frac{P^{(s)}P_{s_k,n}c}{\Vert P^{(s)}P_{s_k,n}c\Vert}\Big\Vert\\&\le&\frac{2\big\Vert P_nc/\Vert P_{s_k}\Vert-P^{(s)}P_{s_k,n}c\big\Vert}{\Vert P_nc\Vert/\Vert P_{s_k}\Vert}\\&\le&\frac{2\varepsilon_k\Vert P_{s_k,n}c\Vert}{\Vert P_nc\Vert/\Vert P_{s_k}\Vert}.\end{array}\end{equation}The r.h.s. tends to 0 in consequence of the hypotheses $c\in\mathcal V_\mathcal A$ and $(H_1)$. To simplify the l.h.s. we prove below that $P^{(s)}$ has rank 1. 

Using the set $J$ and the matrix $Q$ defined in $(H_2)$, $P^{(s)}Q$ is a limit-point of $\frac{P_{s_{k+1}}}{\Vert P_{s_k}\Vert\Vert P_{s_k,s_{k+1}}\Vert}$ and, using $(H_1)$,\begin{equation}\label{PQ}P^{(s)}Q\ne0\text{ and }
\frac{P^{(s)}Q}{\Vert P^{(s)}Q\Vert}=P^{(s)}.\end{equation}

%The convex hull of the normalized columns of $P^{(s)}$ is$$\textstyle\mathcal C:=\{\sum_j\alpha_j\frac{P^{(s)}u_j}{\Vert P^{(s)}u_j\Vert}\;;\;j\text{ such that }P^{(s)}u_j\ne0,\ \alpha_j\ge0,\ \sum_j\alpha_j=1\}.$$
Let $\mathcal C$ be the convex hull of the normalized columns of $P^{(s)}$, that is, the vectors $c_j=\frac{P^{(s)}u_j}{\Vert P^{(s)}u_j\Vert}$ for $j\in J$, and suppose that $\mathcal C$ has two distinct extreme points $v,v'$. From the definition of the extreme points, $v$ cannot be a positive linear combinations of distinct elements of $\mathcal C$, hence there exists $j\in J$ such that $v=c_j$. By (\ref{PQ}), $c_j$ is collinear to $\sum_{i\in J}q_{i,j}c_i$. Since $c_j$ cannot be a positive
linear combination of different vectors $c_i$, and since the $c_i$ have norm 1, $c_j=c_i$ for any $i$ such that $q_{i,j}\ne0$.

Similarly, there exist $j'\in J$ such that $v=c_{j'}$ , and $c_{j'}=c_{i'}$ for any $i'$ such that $q_{i',j'}\ne0$. Using the index $i$ of the hypothesis $(H_2)$, $v=v'=c_i$, a contradiction.

By a theorem of Minkowski (see for instance Theorem 0.4 of 

http://www.mat.unimi.it/users/libor/AnConvessa/ext.pdf), $\mathcal C$ is the convex hull of its extreme points. Hence $\mathcal C$ has only one element $v$, and there exists a row-vector $r$ such that~$P^{(s)}=vr$.

The l.h.s. of (\ref{PnPs}) is $\big\Vert\frac{P_nc}{\Vert P_nc\Vert}-\frac v{\Vert v\Vert}\big\Vert$, so $\frac{P_nc}{\Vert P_nc\Vert}$ converges to $\frac v{\Vert v\Vert}$.\end{proof}

%Theorem \ref{rankone} is made to be used in difficult cases, when the matrices have many zero entries and when no formula gives $P_n$ for any $n\in\mathbb N$. We detail how to use this theorem in the second example of Section \ref{examples}.

%Theorem \ref{rankone} applies obviously if the matrices $A_n$ belong to a finite set $\mathcal M$ of positive matrices: denoting by $\varepsilon$ the smallest entry of the matrices $\frac M{\Vert M\Vert}$, $M\in\mathcal M$, the lower bound in $(H_1)$ is at least $\varepsilon^2$ (when the sequence $(n_k)_{k\in\mathbb N}$ is defined by $n_k=k$), and the $2\times2$ matrices involved in $(H_2)$ are positive. But in this case one can prove also the convergence of $\frac{P_nc}{\Vert P_nc\Vert}$ by using the Birkhoff contraction coefficient \cite[Section 3]{Sen81}, this is why we use rather this method in the first example of Section \ref{examples}.

%In the following examples one can compute $P_n$ and the sequence $\frac{P_nc}{\Vert P_nc\Vert}$ diverges in certain cases; it may converge even if $(H_1)$ or $(H_2)$ does not hold.

\begin{counterexample}Suppose that $A_n\in\{A,B\}$ with $A=\left(\begin{smallmatrix}1&0&1\\0&1&0\\0&0&1\end{smallmatrix}\right)$ and $B=\left(\begin{smallmatrix}1&0&0\\0&1&1\\0&0&1\end{smallmatrix}\right)$. Then $\frac{P_n}{\Vert P_n\Vert}=\left(\begin{smallmatrix}1&0&n_1\\0&1&n_2\\0&0&1\end{smallmatrix}\right)$ where $n_1=n_1(n)$ (resp. $n_2=n_2(n)$) is the number of integers $k\le n$ such that $A_k=A$ (resp. $A_k=B$). Given a positive column-vector $c$, $\frac{P_nc}{\Vert P_nc\Vert}$ diverge if, for instance, $A_n=A\Leftrightarrow\exists k\text{ even }2^k\le n<2^{k+1}$. The necessary and sufficient condition for $\frac{P_nc}{\Vert P_nc\Vert}$ to converge is that $\frac{n_1(n)}{n_2(n)}$ has a finite or infinite limit when $n\to\infty$, and this is also the  necessary and sufficient condition for $(H_1)$ to hold for some limit-point $P^{(s)}$.

More generally the product of $3\times3$ upper-triangular matrices can be easily computed (Appendix D.2) and the necessary and sufficient condition for $\frac{P_nc}{\Vert P_nc\Vert}$ to converge follows if the matrices are nonnegative.\end{counterexample}

\subsection{Approximation by a matrix of rank 1, and a counterexample}\label{approx}We expect that in general, if the matrices $A_n$ are nonnegative there exist a sequence of row-vectors $(r_n)_{n\in\mathbb N}$ and a constant column-vector $v$ such that $\lim_{n\to\infty}\big(\frac{P_n}{\Vert P_n\Vert}-vr_n\big)=0$.

More generally, any sequence of matrices (not necessarily a sequence of product-matrices) satisfies the following straightforward equivalences:

\begin{proposition}\label{matricesMn}(i) For any sequence $(M_n)_{n\in\mathbb N}$ of complex-valued $d\times d$ matrices we have the following equivalences, where $\delta_1(n)\ge\delta_2(n)\ge\dots$ are the singular values of $M_n$:$$\begin{array}{l}\exists v_n\text{ (column-vector)},r_n\text{ (row-vector)}\ \lim_{n\to\infty}\big(\frac{M_n}{\Vert M_n\Vert}-v_nr_n\big)=0\\\Leftrightarrow\lim_{n\to\infty}\frac{\delta_2(n)}{\delta_1(n)}=0\\\Leftrightarrow\text{any limit-point of }\frac{M_n}{\Vert M_n\Vert}\text{ has rank }1.\end{array}$$If it holds, $\lim_{n\to\infty}\det\frac{M_n}{\Vert M_n\Vert}=0$.

(ii) We have the following equivalence if the matrices $M_n$ are nonnegative:$$\begin{array}{l}\text{There exists }v\text{ (constant) and }r_n\text{ such that }\lim_{n\to\infty}\big(\frac{M_n}{\Vert M_n\Vert}-vr_n\big)=0\\\Leftrightarrow\frac{M_nc}{\Vert M_nc\Vert}\text{ converges to a limit independent of }c>0\text{ when }n\to\infty.\end{array}$$The vector $v$ of the l.h.s. is collinear to $\lim_{n\to\infty}\frac{M_nc}{\Vert M_nc\Vert}$ and to the limit of the first left singular vector of $M_n$.\end{proposition}

\begin{proof}[\bf Proof](i)We use the singular value decomposition $M_n=S_n\Delta_n{^t\overline{T_n}}$. If $\lim_{n\to\infty}\frac{\delta_2(n)}{\delta_1(n)}=0$, since $S_n$ and $T_n$ are unitary we deduce the estimation $\frac{M_n}{\Vert M_n\Vert}=v_nr_n+o(1)$ with $v_n=\frac{\delta_1(n)}{\Vert M_n\Vert}S_nu_1$ and $r_n={^tu_1}{^t\overline{T_n}}$. 

On the other hand if $\frac{\delta_2(n)}{\delta_1(n)}$ does not converge to $0$, there exists an increasing sequence $(n_k)_{k\in\mathbb N}$ for the sequences $\frac{\delta_2(n_k)}{\delta_1(n_k)},\frac{M_{n_k}}{\Vert M_{n_k}\Vert},S_{n_k}$ and $T_{n_k}$ to have some limits $\delta$ (nonnull), $M,S$ and $T$ when $k\to\infty$. Hence the diagonal matrix $\frac{\Delta_{n_k}}{\Vert M_{n_k}\Vert}$ has a limit $\Delta$ of rank at least 2. The matrix $M=S\Delta{^t\overline{T}}$ has rank at least 2, hence there cannot exist a column-vector $v_n$ and a row-vector $r_n$ such that $\lim_{n\to\infty}\Big(\frac{M_n}{\Vert M_n\Vert}-v_nr_n\Big)=0$.

This last assertion is equivalent to "$\min_{u\in\mathcal U}\big\Vert\frac{M_n}{\Vert M_n\Vert}-u\big\Vert$ converges to 0 when $n\to\infty$", where $\mathcal U$ is the compact set of the $d\times d$ matrices of norm 1 and rank 1. It is equivalent to "any limit-point of $\frac{M_n}{\Vert M_n\Vert}$ has rank 1".

The equality $\lim_{n\to\infty}\det\frac{M_n}{\Vert M_n\Vert}=0$ follows from $\lim_{n\to\infty}\big(\frac{M_n}{\Vert M_n\Vert}-v_nr_n\big)=0$, by using the uniform continuity of the determinant on the compact set of the $d\times d$ matrices of norm 1.

(ii) Suppose $\exists(v,r_n)\ \lim_{n\to\infty}\big(\frac{M_n}{\Vert M_n\Vert}-vr_n\big)=0$. Multiplying at the right by $\frac{\Vert M_n\Vert}{\Vert M_nc\Vert}c$ (bounded vector), we deduce $\exists \lambda_n\ \lim_{n\to\infty}\big\Vert\frac{M_nc}{\Vert M_nc\Vert}-\lambda_nv\big\Vert=0$. By the second triangular inequality, $\lim_{n\to\infty}\Vert\lambda_n v\Vert=1$, so $\lim_{n\to\infty}\frac{M_nc}{\Vert M_nc\Vert}=\frac v{\Vert v\Vert}$.

Conversely suppose $\exists v=\lim_{n\to\infty}\frac{M_nc}{\Vert M_nc\Vert}$ independent of $c>0$. We apply this to the vector $c=u_{j,\varepsilon}:=u_j+\varepsilon u$, with $\varepsilon>0$:$$\textstyle\exists n_\varepsilon\ \forall n\ge n_\varepsilon\ \big\Vert \frac{M_nu_{j,\varepsilon}}{\Vert M_nu_{j,\varepsilon}\Vert}-c\big\Vert\le\varepsilon.$$Using the classical inequality $\big\Vert\frac x{\Vert x\Vert}-\frac y{\Vert y\Vert}\big\Vert\le\frac{2\Vert x-y\Vert}{\max(\Vert x\Vert,\Vert y\Vert)}$ we deduce$$\textstyle\big\Vert \frac{M_nu_j}{\Vert M_nu_j\Vert}-c\big\Vert\le\frac{2\varepsilon\Vert M_n\Vert}{\Vert M_nu_j\Vert}+\varepsilon\text{ and }\big\Vert \frac{M_nu_j}{\Vert M_n\Vert}-\frac{\Vert M_nu_j\Vert}{\Vert M_n\Vert}c\big\Vert\le3\varepsilon.$$So the l.h.s. of the equivalence (ii) holds with the row-vector $r_n$ of entries $\frac{\Vert M_nu_j\Vert}{\Vert M_n\Vert}$.

To prove the last assertion we use the singular value decomposition with $\lim_{n\to\infty}\frac{\delta_2(n)}{\delta_1(n)}=0$. We deduce that $S_nu_1$ (first left singular vector of $M_n$) is equal to $\frac{M_nc}{\Vert M_nc\Vert_2}+o(1)$, where $\Vert\cdot\Vert_2$ is the Euclidean norm, and consequently converges to $\frac v{\Vert v\Vert_2}$.\end{proof}

\begin{counterexample}In this example (inspired from \cite{Perron}) the sequence $\frac{P_nc}{\Vert P_nc\Vert}$ diverges for any $c\ne0$ and the six assertions, in Proposition \ref{matricesMn}, are false for the sequence $(P_n)_{n\in\mathbb N}$. We consider a sequence of product-matrices $P_n=A_1\cdots A_n$ such that $A_n\in\big\{\left(\begin{smallmatrix}1/2&0\\0&1\end{smallmatrix}\right),\left(\begin{smallmatrix}1&1\\1&0\end{smallmatrix}\right)\big\}$. Equivalently, there exists a sequence of nonnegative integers $(n_k)_{k\in\mathbb N}$ such that\begin{equation}\label{skh}\begin{array}{l}P_n=\left(\begin{smallmatrix}1/2&0\\0&1\end{smallmatrix}\right)^{n_1}\left(\begin{smallmatrix}1&1\\1&0\end{smallmatrix}\right)\cdots\left(\begin{smallmatrix}1/2&0\\0&1\end{smallmatrix}\right)^{n_k}\left(\begin{smallmatrix}1&1\\1&0\end{smallmatrix}\right)\left(\begin{smallmatrix}1/2&0\\0&1\end{smallmatrix}\right)^h\\\text{for any }n=s_k+h\text{ with }s_k=\sum_{i=1}^k(n_i+1)\text{ and }0\le h\le n_{k+1}.\end{array}\end{equation}By induction on $k$, $P_{s_k}=\left(\begin{matrix}v_k&2^{-n_k}v_{k-1}\end{matrix}\right)$ with $v_0=\left(\begin{smallmatrix}1\\0\end{smallmatrix}\right)$, $v_1=\left(\begin{smallmatrix}2^{-n_1}\\1\end{smallmatrix}\right)$ and\begin{equation}\label{vvv}v_k=2^{-n_k}v_{k-1}+2^{-n_{k-1}}v_{k-2}\ (k\ge2).\end{equation}Also true for $k=0$ with, by convention, $P_0=\left(\begin{smallmatrix}1&0\\0&1\end{smallmatrix}\right),n_0=0,v_{-1}=\left(\begin{smallmatrix}0\\1\end{smallmatrix}\right)$.

In consequence of (\ref{vvv}), given a sequence of positive reals $(\epsilon_n)_{n\in\mathbb N}$ we can chose successively the integers $n_k$, $k\ge1$ to have\begin{equation}\label{vk}\textstyle\Vert v_k-2^{-n_{k-1}}v_{k-2}\Vert\le2^{-n_{k-1}}\Vert v_{k-2}\Vert\epsilon_k/2.\end{equation}Using the inequality $\big\Vert\frac x{\Vert x\Vert}-\frac y{\Vert y\Vert}\big\Vert\le\frac{2\Vert x-y\Vert}{\Vert y\Vert}$, the vector $\frac{v_k}{\Vert v_k\Vert}-\frac{v_{k-2}}{\Vert v_{k-2}\Vert}$ has norm at most $\epsilon_k$. Consequently the distance from $\frac{v_k}{\Vert v_k\Vert}$ to one of the two canonical vectors, the first if $n$ is even or the second if $n$ is odd, is at most $\epsilon=\sum_k\epsilon_k$.

Given a nonnull vector $c=\left(\begin{smallmatrix}c_1\\c_2\end{smallmatrix}\right)$, we have also $\Vert P_{s_k}c-c_12^{-n_{k-1}}v_{k-2}\Vert\le\vert c_1+c_2\vert2^{-n_{k-1}}\Vert v_{k-2}\Vert\epsilon_k/2$ and deduce similarly that the distance from $\frac{P_{s_k}c}{\Vert P_{s_k}c\Vert}$ to one of the two canonical vectors, not the same according to the parity of~$k$, is at most $C\epsilon$, where the constant $C$ depends only on the vector $c$. We chose $\epsilon$ small enough and deduce that the sequence $\frac{P_{s_k}c}{\Vert P_{s_k}c\Vert}$ diverges.

Using (\ref{skh}), (\ref{vvv}) and (\ref{vk}), the first column of $P_{s_k}$ (resp. of $P_{s_k+n_{k+1}}$) has norm larger (resp. smaller, if $\epsilon_k<2$) than the second. Hence there exists an integer $h_k$ for the norms of the columns of $P_{s_k+h_k}$to have ratio at least 1/2 and at most~2. Consequently, given $\alpha>0$ we can chose $\epsilon$ for the diagonal entries of $\frac{P_{s_k+h_k}}{\Vert P_{s_k+h_k}\Vert}$ ($k$ even) to belong to $[\frac13-\alpha,\frac23+\alpha]$, and the other entries to $[0,\alpha]$. Choosing $\alpha$ small enough, the determinant of $\frac{P_{s_k+h_k}}{\Vert P_{s_k+h_k}\Vert}$ is bounded from 0 and the six assertions in Proposition~\ref{matricesMn} are false for the sequence $(P_n)_{n\in\mathbb N}$.\end{counterexample}

\begin{remark}About the singular values of a product of triangular matrices: let $P_n=\prod_{k=1}^n\left(\begin{smallmatrix}a_k&b_k\\0&d_k\end{smallmatrix}\right)$ with $\forall k\ a_kd_k\ne0$. Since $P_n=\left(\begin{smallmatrix}a^*_n&d^*_ns_n\\0&d^*_n\end{smallmatrix}\right)$ with $a^*_n=a_1\cdots a_n$, $d^*_n=d_1\cdots d_n$, $s_n=\sum_1^n\frac{a^*_k}{d^*_k}\frac{b_k}{a_k}$, its singular values represent the orders of growth of its columns iff $\vert s_n\vert$ is bounded, because the first singular value has the same order of growth as the modulus of the largest entry, and their product is the modulus of $\det P_n$.

Let $M_n=(m_{i,j}^{(n)})_{1\le i\le d\atop1\le j\le d}$ be a sequence of $d\times d$ upper-triangular nonnegative matrices with positive diagonal entries. We prove that each singular value $\delta_k(n)$ has the same order of growth as the $k^{\rm th}$ columns of $M_n$ iff\begin{equation}\label{sin}\textstyle\sup\big\{\frac{m_{i,j}^{(n)}}{m_{i',j'}^{(n)}}\;;\;i,i',j,j',n\text{ such that }j\le j'\big\}<\infty.\end{equation}For this we consider the characteristic polynomial of ${^tM_n}\cdot M_n$$$P(X)=\textstyle\det(XI-{^tM_n}\cdot M_n)=\prod_k(X-\delta_k(n)^2).$$For any $k$, $(\delta_1(n)\cdots\delta_k(n))^2$ has the same order of growth as the coefficient $a_k$ of $X^{d-k}$ in $P(X)$, because $\delta_1(n)\ge\delta_2(n)\ge\dots$. And $(-1)^ka_k$ is the sum of the positive determinants $\det{^tM_I}\cdot M_I$, where $M_I$ is the $\{1,\dots,d\}\times I$ submatrix of $M_n$, for any $I\subset\{1,\dots,d\}$ of cardinality $k$. It is the following sum for $I\subset\{1,\dots,d\}$, $\#I=k,\ \sigma$ permutation 
of $I$, $(j_1,\dots,j_k)\in I^k$:$$(-1)^ka_k=\sum_{I,\sigma,j_1,\dots, j_k}(-1)^{s(\sigma)}m_{j_1,i_1}m_{j_1,\sigma(i_1)}\cdots m_{j_k,i_k}m_{j_k,\sigma(i_k)}.$$One can assume that $j_1,\dots,j_k$ are distinct because any term such that $\exists\alpha,\beta\ j_\alpha=j_\beta$, is the opposite of the term corresponding to the permutation $\sigma\circ\tau$, $\tau$ transposition of $i_\alpha$ and $i_\beta$. So if (\ref{sin}) is true, this sum belongs to $O((m_{1,1}\cdots m_{k,k})^2)$. 

Since the $\{1,\dots,k\}\times\{1,\dots,k\}$-submatrix of ${^tM_n}\cdot M_n$
is ${^tM}\cdot M$ (where $M$ is the $\{1,\dots,k\}\times\{1,\dots,k\}$-submatrix of $M_n$) and has determinant $(m_{1,1}\cdots m_{k,k})^2$, we deduce $\forall k\ (\delta_1(n)\cdots\delta_k(n))^2\asymp(m_{1,1}\cdots m_{k,k})^2$ hence $\forall k\ \delta_k(n)\asymp m_{k,k}\asymp\Vert M_nu_k\Vert$.

Conversely if the ratios $\frac{m_{i,j}^{(n)}}{m_{i',j'}^{(n)}}$ are not bounded when $j\ge j'$, either there exists $j$ for $\frac{\Vert M_nu_j\Vert}{m_{j,j}}$ to be not bounded, or there exists two indices $j\ge j'$ for $\frac{\Vert M_nu_j\Vert}{\Vert M_nu_{j'}\Vert}$ to be not bounded. In the first case the relation $\delta_1(n)\cdots\delta_d(n)\asymp\Vert M_nu_1\Vert\cdots \Vert M_nu_d\Vert$ cannot hold, and in the second, since $\delta_j(n)\le\delta_{j'}(n)$, the relation $\delta_i(n)\asymp\Vert M_nu_i\Vert$ is false for $i=j$ or for $i=j'$.\end{remark}

\section{Applications of the infinite products of matrices}\label{functional}

\subsection{Two-scale difference equations}\label{sltsde}The "two-scale difference equations" \cite{DL91,DL92,DL92Two,BeP} are equations in $f:\mathbb R\to\mathbb R$, of the form\begin{equation}\label{ltsde}\textstyle f(x)=\sum_{n=0}^Nc_nf(kx-n)\text{ for any }x\in\mathbb R,\text{ where }c_n\in\mathbb R\text{ and }1<k\in\mathbb R\end{equation} (lattice two scale if $k\in\mathbb N$). If $k\in\mathbb N$ one usually solve this equation by means of the fundamental curves~\cite{MP,BW}. The fundamental curve associated to a set of matrices $\{B_0,\dots,B_{k-1}\}$ is parametrized by a function $\psi:[0,1]\to\mathbb R^d$ that satisfies, for any $0\le i\le k-1$ and any $x$ in the closed interval $[\frac ik,\frac{i+1}k]$,\begin{equation}\label{fcopenorclose}\textstyle\psi(x)=B_i\psi(kx-i).\end{equation}More precisely one can solve (\ref{ltsde}) by applying both following propositions:

\begin{proposition}\label{fcltsde}Let $q:=\left\lceil\frac N{k-1}\right\rceil-1$ and $B_i:=\left(\begin{smallmatrix}c_i&c_{i-1}&\dots&c_{i-q}\\c_{k+i}&c_{k+i-1}&\dots&c_{k+i-q}\\\vdots&\vdots&\ddots&\vdots\\c_{kq+i}&c_{kq+i-1}&\dots&c_{kq+i-q}&\end{smallmatrix}\right)$ with $0\le i\le k-1$ and $c_0,\dots,c_N\in\mathbb R$, $c_n=0$ for $n\not\in\{0,\dots,N\}$. A function $f:\mathbb R\to\mathbb R$, null outside $[0,q+1)$, is solution of (\ref{ltsde}) iff the function $\psi_f$ defined on $[0,1)$ by $\psi_f(x):=\left(\begin{smallmatrix}f(x)\\f(x+1)\\\vdots\\f(x+q)\end{smallmatrix}\right)$ satisfies (\ref{fcopenorclose}) for any $0\le i\le k-1$ and any $x$ in the semi-open interval $[\frac ik,\frac{i+1}k)$.\end{proposition}\begin{proof}[\bf Proof]A function $f:\mathbb R\to\mathbb R$ null outside $[0,q+1)$ is solution of (\ref{ltsde}) iff, for any $x\in[0,1)$ and for any $j\in\{0,\dots,q\}$, the $(j+1)^{\rm th}$ entry of $\psi_f(x)$ is equal to $\sum_{n\in\mathbb Z}c_nf(k(x+j)-n)$. For any $x\in[0,1)$ there exist $x'\in[0,1)$ and $0\le i\le k-1$ such that $x=\frac{x'+i}k$. By change of variable, one has$$\textstyle\sum_{n\in\mathbb Z}c_nf(k(x+j)-n)=\sum_{j'\in\mathbb Z}c_{i+kj-j'}f(x'+j').$$The r.h.s. is the $(j+1)^{\rm th}$ entry of $B_i\psi_f(x')$ because $f(x'+j')$ is null  when $j'\not\in\{0,\dots,q\}$. Hence one has the equivalence: $f\text{ solution of (\ref{ltsde}) }\Leftrightarrow\psi_f(x)=B_i\psi_f(x')=B_i\psi_f(kx-i)$.\end{proof}

The second proposition is straightforward:\begin{proposition}\label{openorclose}Suppose that, for any $x=\sum_{n=1}^\infty\frac{x_n}{k^n}$ (expansion in base $k$), there exists $\psi(x)=\lim_{n\to\infty}B_{x_1}\cdots B_{x_n}c$, where $c$ is a $d$-dimensional real column-vector. Then (\ref{fcopenorclose}) holds for any $x$ in the semi-open interval $[\frac ik,\frac{i+1}k)$, $0\le i\le ,k-1$.\end{proposition}See \cite[Theorem 2.1]{MP} or \cite[Theorem 6.1]{DL92Two} for examples of such functions.

\begin{remark}Proposition \ref{fcltsde} considers only the solutions null outside $[0,q+1)$ because, from \cite[Corollary~2.2]{DL91}, any integrable and continuous solution of~(\ref{ltsde}) is null outside this interval. On the other hand, the condition on $\psi_f$ in Proposition \ref{fcltsde} does not guarantee that the curve parametrized by $\psi_f$ is a fundamental curve; it is a fundamental curve iff (\ref{fcopenorclose}) holds for any $x$ in the closed interval $[\frac ik,\frac{i+1}k]$, that is, iff for any $i$ $\lim_{n\to\infty}B_i(B_{k-1})^nc=\lim_{n\to\infty}B_{i+1}(B_0)^nc$.\end{remark}

\subsection{Case where the equation (\ref{ltsde}) has a nonnegative solution in $L^1$}\label{L1}By integrating (\ref{ltsde}) from $-\infty$ to $x$ one obtains, setting $p_n:\frac {c_n}k$,\begin{equation}\label{iltsde}\textstyle F(x)=\sum_{n=0}^Np_nF(kx-n).\end{equation}So, searching a solution of (\ref{ltsde}) with the conditions "$0\le f\in L^1$ and $f$ not a.e. null", is equivalent to search a solution of (\ref{iltsde}) with the conditions "$F$ absolutely continuous and non-decreasing from 0 to a positive real".

One can also consider (\ref{iltsde}) as an equation in $F$, with only the conditions "$F$ left-continuous (resp. right-continuous) and non-decreasing from 0 to a positive real". If $F$ exists one has necessarily $\sum_np_n=1$.

According to the introduction of \cite{H}, if $p_n>0$ and $\sum_np_n=1$ this equation has a unique solution, up to a multiplicative constant,  because it is equivalent to search a finite Borel measure $\mu$ such that, for any Borel set~$B$,\begin{equation}\label{integrated}\textstyle\mu(B)=\sum_{n=0}^Np_n\mu(kB-n).\end{equation}.

\begin{remark}\label{notac}According to \cite[Proposition 5.3]{JMNA}, if $p_0=\dots=p_N$ and \hbox{$k\in\mathbb N$}, the necessary and sufficient condition for $\mu$ to be absolutely continuous with respect to the Lebesgue measure, is that $k$ divides $N+1$. So, if $k$ does not divide $N+1$, (\ref{ltsde}) does not have a solution $0\le f\in L^1$ except the a.e. null ones. In such a case, from Proposition \ref{openorclose} the limit of $B_{x_1}\cdots B_{x_n}c$, if it exists for any $x$ of expansion $\sum_{i=1}^\infty\frac{x_n}{k^n}$, is a.e. null. However (\ref{ltsde}) has a.e. null solutions, for instance the function $f$ defined by $f(x)=0$ for $x\not\in[0,q)$ and $\psi_f(x)=\left\{\begin{array}{ll}B_{x_1}\cdots B_{x_n}c&\text{when }x\text{ has a finite expansion }\sum_{i=1}^n\frac{x_i}{k^i}\\0&\text{otherwise}\end{array}\right.$, where $c$ is a column-vector.\end{remark}

We are interested by the case where $\mu$ is not absolutely continuous. In this case we do not apply Proposition~\ref{fcltsde} to the derivative of the distribution function of~$\mu$, because from \cite{JW} this derivative is a.e. null, but we use the following similar proposition:\begin{proposition}\label{Fcltsde}Let $q:=\left\lceil\frac N{k-1}\right\rceil-1$ and $M_i:=\left(\begin{smallmatrix}p_i&p_{i-1}&\dots&p_{i-q}\\p_{k+i}&p_{k+i-1}&\dots&p_{k+i-q}\\\vdots&\vdots&\ddots&\vdots\\p_{kq+i}&p_{kq+i-1}&\dots&p_{kq+i-q}&\end{smallmatrix}\right)$, $0\le i\le k-1$ ($k\in\mathbb N$), where the coefficients $p_0,\dots,p_N$ are reals and  $p_n=0$ for any $n\not\in\{0,\dots,N\}$. Any solution $F$ of (\ref{iltsde}) which is constant on $(-\infty,0)$ and constant on $[q+1,+\infty)$ satisfies the following relation, when the reals $x$ and $y$ are in a same interval $[\frac ik,\frac{i+1}k)$, $0\le i\le k-1$:$$\textstyle\phi_F(x)-\phi_F(y)=M_i\big(\phi_F(kx-i)-\phi_F(ky-i)\big).$$\end{proposition}\begin{proof}[\bf Proof]For $x=\frac{x'+i}k$ and $y=\frac{y'+i}k$ with $x',y'\in[0,1)$ and $0\le i\le k-1$, the $(j+1)^{\rm th}$ entry of $\phi_F(x)-\phi_F(y)$ is $\sum_{n\in\mathbb Z}p_n\big(F(k(x+j)-n)-F(k(y+j)-n)\big)$. By change of variable it is $\sum_{j'\in\mathbb Z}p_{i+kj-j'}\big(F(x'+j')-F(y'+j')\big)$, that is, the $(j+1)^{\rm th}$ entry of $M_i\psi_F(x')-M_i\psi_F(y')$ because $F(x'+j')-F(y'+j')$ is null  when $j'\not\in\{0,\dots,q\}$.\end{proof}

So if $\mu$ is a solution of (\ref{integrated}) then, for any interval of the form $I=I_{x_1,\dots,x_n}:=\big[\sum_1^n\frac{x_i}{k^i},\sum_1^n\frac{x_i}{k^i}+\frac1{k^n}\big)$ (set of the reals whose expansion in base $k$ begins by $x_1,\dots,x_n$) one has\begin{equation}\label{linearep}\textstyle\left(\begin{smallmatrix}\mu(I)\\\mu(I+1)\\\vdots\\\mu(I+q)\end{smallmatrix}\right)=\prod_1^n M_{x_i}c\text{ with }c:=\left(\begin{smallmatrix}\nu([0,1))\\\nu([1,2))\\\vdots\\\nu([q,q+1))\end{smallmatrix}\right).\end{equation}One can compute $c$ as an eigenvector of $\sum_iM_i$ (consequence of the formula in the case $n=1$). This formula allows to study the weak-Gibbs property of $\mu$ in the first example of Section \ref{examples}.

\subsection{Bernoulli convolutions}\label{Bc}By generalization of \cite{Erd39,PSS} we call "the Bernoulli convolution associated to the real $k>1$ and the probability vector $p=(p_0,\dots,p_N)$", the probability measure $\nu_{k,p}$ solution of (\ref{integrated}). According to the introduction of \cite{H}, the support of $\nu_{k,p}$ is a subset of $\big[0,\frac N{k-1}\big]$.

We give an equivalent definition: let $\nu_{k,p}(I)$ ($I$ interval) be the probability that $\sum_{n=1}^\infty\frac{\omega_n}{k^n}\in I$, when the $\omega_n$ take each value $i\in\{0,\dots,N\}$ with probability $p_i$ and are independent; it is clear that $\nu_{k,p}$ satisfies (\ref{integrated}).

To give an example with $k\not\in\mathbb N$, we consider the Bernoulli convolution associated to the real $\beta>1$ such that $\beta^3=2\beta^2-\beta+1$ and to the probability vector $p=(\frac12,\frac12)$.
A matrix-representation of $\nu_{\beta,2}$ is obtained as follows:

\begin{lemma}\label{SiM}We consider the three words $w_0:=0$, $w_1:=10$, $w_2:=1100$. For any $i_1\dots i_n\in\{0,1,2\}^n$, let $\varepsilon_1\dots\varepsilon_m\in\{0,1\}^m$ be the concatenation of the words $w_{i_1},\dots, w_{i_n}$. We denote by $I^{(\beta)}_{i_1,\dots,i_n}$ the interval $\big[\sum_{i=1}^m\frac{\varepsilon_i}{\beta^i},\frac1{\beta^m}+\sum_{i=1}^m\frac{\varepsilon_i}{\beta^i}\big)$, set of the reals whose expansion in base $\beta$ begins by $\varepsilon_1\dots\varepsilon_m$. The matrix-representation of $\nu_{\beta,2}$ is\begin{equation}\label{matre}\textstyle\nu_{\beta,2}\big(\frac1{\beta-1}I^{(\beta)}_{i_1,\dots,i_n}\big)=\left\{\begin{array}{ll}{^tu_1}M'_{i_2}\cdots M'_{i_n}c_\beta&\text{if }i_1=0\\\frac12\ {^tu_3}M'_{i_2}\cdots M'_{i_n}c_\beta&\text{if }i_1=1\\\frac18\ {^tu_5}M'_{i_2}\cdots M'_{i_n}c_\beta&\text{if }i_1=2\end{array}\right.\end{equation}with

$\scriptstyle2M'_0:=\left(\begin{smallmatrix}1&0&0&0&0&0&0\\0&0&1&0&0&0&0\\0&0&0&1&1&0&0\\0&0&0&0&0&0&0\\1&0&0&0&0&0&1\\0&0&0&0&1&0&0\\0&1&0&0&0&0&0\end{smallmatrix}\right),\ 4M'_1:=\left(\begin{smallmatrix}0&0&1&1&0&0&0\\0&0&0&0&0&1&0\\0&0&0&1&1&0&0\\1&0&0&0&0&0&0\\0&0&1&0&0&0&0\\0&0&0&0&0&0&0\\0&0&0&0&0&0&0\end{smallmatrix}\right),\ 16M'_2:=\left(\begin{smallmatrix}1&0&0&0&1&0&1\\0&0&0&0&0&0&0\\1&0&0&0&0&0&1\\0&0&0&1&1&0&0\\0&0&0&0&1&0&0\\0&0&0&0&0&0&0\\0&0&0&0&0&0&0\end{smallmatrix}\right),\ 20c_\beta:=\left(\begin{smallmatrix}12\\8\\13\\4\\12\\6\\4\end{smallmatrix}\right)$.\end{lemma}

\begin{proof}[\bf Proof]Let $j_1:=0,j_2:=\beta-1,j_3:=\frac1{\beta^2},j_4:=\frac{1-\beta}\beta,j_5:=\frac1\beta,j_6:=\frac{\beta+1}{\beta^2},$

$j_7=1$, we first prove that\begin{equation}\label{translatedcylinders}\left(\begin{smallmatrix}\nu_{\beta,2}\big(j_1+(\beta-1)I^{(\beta)}_{i_1,\dots,i_n}\big)\\\vdots\\\nu_{\beta,2}\big(j_7+(\beta-1)I^{(\beta)}_{i_1,\dots,i_n}\big)\end{smallmatrix}\right)=M'_{i_1}\cdots M'_{i_n}\left(\begin{smallmatrix}\nu_{\beta,2}\big(j_1+(\beta-1)[0,1)\big)\\\vdots\\\nu_{\beta,2}\big(j_7+(\beta-1)[0,1)\big)\end{smallmatrix}\right).\end{equation}Suppose for instance that $i_1=1$ (the cases $i_1=0$ and $i_1=2$ are similar). Then $I^{(\beta)}_{i_1,\dots,i_n}=\frac1\beta+\frac1{\beta^2}I^{(\beta)}_{i_2,\dots,i_n}$. The random variable $\sum_{n=1}^\infty\frac{\omega_n}{\beta^n}$ belongs to $j+(\beta-1)I^{(\beta)}_{i_1,\dots,i_n}$ ($j\in\mathbb R$) iff$$\textstyle\sum_{n=3}^\infty\frac{\omega_n}{\beta^{n-2}}\in j'+(\beta-1)I^{(\beta)}_{i_2,\dots,i_n}\ \text{\footnotesize with }j'=-\omega_1\beta-\omega_2+\beta^2j+(\beta-1)\beta.$$$\omega_1$, $\omega_2$ take the values 0 and 1 with probability $\frac12$, and $j'=j'(\omega_1,\omega_2,j)$ belongs to the interval $\big(-(\beta-1),\frac1{\beta-1})$ because $\sum_{n=3}^\infty\frac{\omega_n}{\beta^{n-2}}\in\big(0,\frac1{\beta-1})$.

We prove now that if $\exists i\ j=j_i$ then $\exists i'\ j'=j_{i'}$: indeed for $i=1$ one has $i'\in\{3,4\}$; and for $i=2$ one has $i'=6$; and so one for $i=3,\dots,7$. So the couples $(i,i')$ are the indices of the nonnull entries in $M'_1$.

Consequently, if the formula (\ref{translatedcylinders}) is true for $I^{(\beta)}_{i_2,\dots,i_n}$, it is true for $I^{(\beta)}_{i_1,\dots,i_n}$. Since (\ref{translatedcylinders}) is obvious when $\i_1\dots\i_n$ is the empty word, it is true for any $n$.

The column-vector at the right of (\ref{translatedcylinders}) is an eigenvector of $M'_0+M'_1+M'_2$ (by applying this formula to $n=1$). To prove that this vector is $c_\beta$, it is sufficient to note that the sum of its first two entries is $\nu_\beta([0,2(\beta-1)))=1$.

(\ref{matre}) follows from (\ref{translatedcylinders}) because$$\textstyle\frac1{\beta-1}I^{(\beta)}_{i_1,\dots,i_n}=\left\{\begin{array}{lcl}j_1+(\beta-1)I^{(\beta)}_{i_2,\dots,i_n}&\text{if }\i_1=0\\j_2+(\beta-1)I^{(\beta)}_{0,i_2,\dots,i_n}&\text{if }\i_1=1\\j_2+(\beta-1)I^{(\beta)}_{1,0,i_2,\dots,i_n}&\text{if }\i_1=2\end{array}\right.$$\end{proof}

\section{Convergence of the sequence $\frac{P_nc}{\Vert P_nc\Vert}$ and weak-Gibbs property in two examples}\label{examples}

\subsection{First example}\label{firstexample}

We say that a measure $\nu$ on the interval $[0,1]$ has the week-Gibbs property with respect to the numeration in the integral base $k$ if there exists a continuous map $\Phi:[0,1]\to\mathbb R$, called a potential of $\nu$, such that
\begin{equation}\label{defGibbs}
\lim_{n\to\infty}\Big(\frac{\nu(I_n(x))}{\exp\big(\sum_{i=0}^{n-1}\Phi(T^i(x)\big)}\Big)^{\frac1n}=1\end{equation}uniformly on $x\in[0,1)$, where $I_n(x)$ denotes the set of the reals $y\in[0,1)$ of same digit as $x$ until the rank $n$ ($\forall i\le n\ y_i=x_i$).

Equivalently, $\nu$ has the weak-Gibbs property in the sense of Definition~\ref{2definitions}(ii) below, with respect to the transformation $T$ defined by $T(x):=\sum_{n=1}^\infty\frac{x_{n+1}}{k^n}$ for any $x$ of representation $x=\sum_{n=1}^\infty\frac{x_n}{k^n}$ in base $k$
($T(x)=kx-\lfloor kx\rfloor$).

This notion depends on the choice of the base, but the aim is to deduce that the multifractal formalism holds for the measure $\nu$, and this property depends only on the measure, see Appendix \ref{mf} for more details.

We consider the measure $\nu_{k,p}$ defined in Subsection \ref{Bc} when $k\ge2$ is an integer, the the matrices $M_i$ defined in Proposition~\ref{Fcltsde}, and a positive vector~$c$. There exist two integers $q$ and $1\le r\le k-1$ such that $N=(k-1)q+r$, and we have\begin{equation}\label{0ri}M_0=:\left(\begin{smallmatrix}p_0&0\\v_0&M'_0\end{smallmatrix}\right),\ M_r=:\left(\begin{smallmatrix}M'_r&v_r\\0&p_N\end{smallmatrix}\right),\ M_i=:\left(\begin{smallmatrix}M'_i&v_i\\0&0\end{smallmatrix}\right)\ (i>r).\end{equation}Since $N\ge k$ (otherwise the matrices $M_i$ have order 1), the vector $v_0=\left(\begin{smallmatrix}p_k\\\vdots\\p_{kq}\end{smallmatrix}\right)$ is nonnull.

\begin{theorem}\label{cvandGibbs}(i) Under the following hypothesis:$$\textstyle(H):p_N\le\rho(M'_r)\text{ (spectral radius of }M'_r\text{) or }\frac N{k-1}\in\mathbb N,$$the sequence $c_{\omega,n}=\frac{M_{\omega_1}\cdots M_{\omega_n}c}{\Vert M_{\omega_1}\cdots M_{\omega_n}c\Vert}$ converges uniformly on $\omega\in\{0,\dots,k-1\}^{\mathbb N}$.

(ii) For any $h\in\{0,\dots,q\}$ the measure $\nu_{k,p}(\cdot+h)$ restricted to the interval $[0,1]$ has the weak-Gibbs property with respect to the numeration in base $k$ if and only if $(H)$ holds.\end{theorem}

The proof is given in Appendix \ref{proofcvandGibbs}.

\begin{remark}If $(H)$ does not hold, the sequence $\frac{M_{\omega_1}\cdots M_{\omega_n}c}{\Vert M_{\omega_1}\cdots M_{\omega_n}c\Vert}$ may converge pointwise but not uniformly. In the special case $N=k=3$, $M_0=\left(\begin{smallmatrix}p_0&0\\p_3&p_2\end{smallmatrix}\right)$, $M_1=\left(\begin{smallmatrix}p_1&p_0\\0&p_3\end{smallmatrix}\right)$, $M_2=\left(\begin{smallmatrix}p_2&p_1\\0&0\end{smallmatrix}\right)$, if $p_3>p_1$, the non-uniformity of the convergence is a consequence of the discontinuity of the limit-function on the set $\{0,\dots,k-1\}^{\mathbb N}$: indeed $\lim_{n\to\infty}\frac{{M_1}^nc}{\Vert{M_1}^nc\Vert}=\left(\begin{smallmatrix}\frac1{1+p_3-p_1}\\\frac{p_3-p_1}{1+p_3-p_1}\end{smallmatrix}\right)$ is not the limit when $m\to\infty$ of $\lim_{n\to\infty}\frac{{M_1}^m{M_2}^{n-m}c}{\Vert{M_1}^m{M_2}^{n-m}c\Vert}=\left(\begin{smallmatrix}1\\0\end{smallmatrix}\right)$.\end{remark}

\subsection{Second example}\label{appli}

Let $\beta$ be the real solution of the equation $\beta^3=2\beta^2-\beta+1$ and let $T:x\mapsto\left\{\begin{array}{ll}\beta x&\text{if }x<\frac1\beta\\\beta^2 (x-\frac1\beta)&\text{if }\frac1\beta\le x<\frac1\beta+\frac1{\beta^2}\\\beta^4(x-\frac1\beta-\frac1{\beta^2})&\text{if }\frac1\beta+\frac1{\beta^2}\le x<1.\end{array}\right.$\begin{theorem}\label{betareal}Given a positive 7-dimensional column-vector $c$ and$$M_0:=\left(\begin{smallmatrix}1&0&0&0&0&0&0\\0&0&1&0&0&0&0\\0&0&0&1&1&0&0\\0&0&0&0&0&0&0\\1&0&0&0&0&0&1\\0&0&0&0&1&0&0\\0&1&0&0&0&0&0\end{smallmatrix}\right),\ M_1:=\left(\begin{smallmatrix}0&0&1&1&0&0&0\\0&0&0&0&0&1&0\\0&0&0&1&1&0&0\\1&0&0&0&0&0&0\\0&0&1&0&0&0&0\\0&0&0&0&0&0&0\\0&0&0&0&0&0&0\end{smallmatrix}\right),\ M_2:=\left(\begin{smallmatrix}1&0&0&0&1&0&1\\0&0&0&0&0&0&0\\1&0&0&0&0&0&1\\0&0&0&1&1&0&0\\0&0&0&0&1&0&0\\0&0&0&0&0&0&0\\0&0&0&0&0&0&0\end{smallmatrix}\right),$$(i) the sequence $c_{\omega,n}=\frac{M_{\omega_1}\cdots M_{\omega_n}c}{\Vert M_{\omega_1}\cdots M_{\omega_n}c\Vert}$ converges uniformly on $\{0,1,2\}^{\mathbb N}$.% and the measure $\nu\big(\frac \cdot{\beta-1}\big)$ has the weak-Gibbs property with respect to the numeration in base $\beta$;

(ii) the measure $\nu_{\beta,2}\big(\frac \cdot{\beta-1}\big)$ has the weak-Gibbs property with respect to $T$ (in the sense of Definition~\ref{2definitions}(ii)).
\end{theorem}

The proof is given in Appendix \ref{proofbetareal}.

\section{Divergence of the sequence of matrices $\frac{P_n}{\Vert P_n\Vert}$}\label{limpoints}

We first point out that, if the matrices $A_n$ are positive and if the sequence $\frac{A_n}{\Vert A_n\Vert}$ converges to a positive matrix, then the sequence $\frac{P_n}{\Vert P_n\Vert}$ converges to
a rank-one matrix. For this we apply \cite[Theorem~1.1]{F} to the matrices $B_n=\frac{^t{A_n}}{\Vert A_n\Vert}$: the sequence $\frac{B_n\cdots B_1}{\Vert B_n\cdots B_1\Vert}$ converges to a rank-one matrix, hence by transposition the sequence $\frac{A_1\cdots A_n}{\Vert ^t{A_n}\cdots^t{A_1}\Vert}$ converges to a rank-one matrix. We deduce that $\frac{\Vert A_1\cdots A_n\Vert}{\Vert ^t{A_n}\cdots^t{A_1}\Vert}$ converges to a positive number hence $\frac{A_1\cdots A_n}{\Vert A_1\cdots A_n\Vert}$ converges to a rank-one matrix.

Theorem~\ref{aediv} below can be seen as a partial converse of this result: for any sequence $(A_n)_{n\in\mathbb N}$ of complex-valued matrices such that $\frac{A_1\cdots A_n}{\Vert A_1\cdots A_n\Vert}$ converges, either the sequence $\frac{A_n}{\Vert A_n\Vert}$ converges, or this sequence has several limit-points with a common left-eigenvector (the proof is similar to the one of Elsner and Friedland in \cite[Theorem~1]{EF97}).

\begin{theorem}\label{aediv}Let $(A_n)_{n\in\mathbb N}$ be a sequence of $d\times d$ matrices on $\mathbb C$, and suppose that the product-matrix $P_n=A_1\cdots A_n$ is nonnull for any $n$. The sequence $\frac{P_n}{\Vert P_n\Vert}$ diverges if the sequence $\frac{A_n}{\Vert A_n\Vert}$ diverges and if its limit-points do not have a common left-eigenvector.\end{theorem}\begin{proof}[\bf Proof]We prove this assertion by contraposition. We first prove that:\begin{lemma}\label{aedivlem}If the sequence $\frac{P_n}{\Vert P_n\Vert}$ converges, there exist a row-vector $r\ne0$ and a sequence of positive reals $(\lambda_n)_{n\in\mathbb N}$ such that\begin{equation}\label{ralambda}\textstyle\lim_{n\to\infty}r\big(\frac{A_n}{\Vert A_n\Vert}-\lambda_nI_d\big)=0.\end{equation}\end{lemma}\begin{proof}[\bf Proof] (i) Let $\lambda_n:=\frac{\Vert P_n\Vert}{\Vert P_{n-1}\Vert\cdot\Vert A_n\Vert}$. The bounded sequence $\lambda_n$ satisfies$$\textstyle\lambda_n\Big(\frac{P_n}{\Vert P_n\Vert}-\frac{P_{n-1}}{\Vert P_{n-1}\Vert}\Big)=\frac{P_{n-1}}{\Vert P_{n-1}\Vert}\Big(\frac{A_n}{\Vert A_n\Vert}-\lambda_nI_d\Big)$$hence, if the sequence $\frac{P_n}{\Vert P_n\Vert}$ converges to a matrix $P$, $\lim_{n\to\infty}\frac{P_{n-1}}{\Vert P_{n-1}\Vert}\big(\frac{A_n}{\Vert A_n\Vert}-\lambda_nI_d\big)=0$ and $\lim_{n\to\infty}P\big(\frac{A_n}{\Vert A_n\Vert}-\lambda_nI_d\big)=0$. Since $P$ has norm $1$, it has at least one non-null row $r$ and (\ref{ralambda}) holds.\end{proof}If the sequence $\frac{P_n}{\Vert P_n\Vert}$ converges the relation (\ref{ralambda}) implies, by the second triangular inequality, $\lim_{n\to\infty}\big(\big\Vert r\frac{A_n}{\Vert A_n\Vert}\big\Vert-\lambda_n\Vert r\Vert\big)=0$. So for any limit-point $A=\lim_{k\to\infty}\frac{A_{n_k}}{\Vert A_{n_k}\Vert}$, the sequence $k\mapsto\lambda_{n_k}$ converges to $\lambda=\frac{\Vert rA\Vert}{\Vert r\Vert}$ and, by (\ref{ralambda}), $r$ is a left-eigenvector of $A$. Either the sequence $\frac{A_n}{\Vert A_n\Vert}$ converges, or it has several limit-points with $r$ as common left-eigenvector.\end{proof}

To prove that the sequence $\frac{P_n}{\Vert P_n\Vert}$ diverges with probability 1, the probability we use can be any infinite product of a non-atomic Borel probability measure $p$ with support $\mathbb C$. This means that the probability for the $(i,j)^{\rm th}$ entry of $A_n$ to belong to the Borel set $B_{n,i,j}$ for any $n,i,j$, is $\prod_{n,i,j}p(B_{n,i,j})$.\begin{corollary}Let $p$ be a non-atomic Borel probability measure with support $\mathbb C$, and $p^*$ the product-probability on $\mathcal M_d(\mathbb C)^\mathbb N$. For $p^*$-a.e. sequence $(A_n)_{n\in\mathbb N}\in\mathcal M_d(\mathbb C)^\mathbb N$, the matrices $P_n$ are nonnull and the sequence $\frac{P_n}{\Vert P_n\Vert}$ diverges.\end{corollary}\begin{proof}[\bf Proof]We have to prove that the following sets have probability~$0$:
$$
\begin{array}{rcl}
\mathcal E_0&:=&\{(A_n)_{n\in\mathbb N}\;;\;\exists n,\ P_n=0\},\\
\mathcal E&:=&\big\{(A_n)_{n\in\mathbb N}\;;\;\text{the limit-points of }\frac{A_n}{\Vert A_n\Vert}\text{ have at least one common}\\&&\text{left-eigenvector}\big\}.\end{array}
$$$\mathcal E_0$ is a subset of $\{(A_n)_{n\in\mathbb N}\;;\;\exists n\ \det A_n=0\}$. It is sufficient to prove that the set $\mathcal E_d:=\{M\in\mathcal M_d(\mathbb C)\;;\;\det M=0\}$ has probability $0$. Suppose by induction hypothesis that $\mathcal E_{d-1}$ has probability $0$. This implies that the co-matrix $M'$ of a matrix $M\in\mathcal M_d(\mathbb C)$ has non-null $(1,1)$-entry with probability~$1$. If this entry is non-null and if $\det M=0$, the $(1,1)$-entry of $M$ is a function of the other entries of $M$:\begin{equation}\label{m11mdd}m_{1,1}=f(m_{1,2},\dots,m_{d,d}).\end{equation}Since $p$ is non-atomic, the conditional probability that (\ref{m11mdd}) holds, given the entries $m_{i,j}$ for $(i,j)\ne(1,1)$, is null. So $\mathcal E_d$ has probability $0$ and $p^*(\mathcal E_0)=0$.

To prove that $p^*(\mathcal E)=0$ we chose two normalized matrices without common left-eigenvector:$$T^{(0)}:=\frac2{d(d+1)}\left(\begin{smallmatrix}1&0&\dots&0\\1&1&\dots&0\\\vdots&\vdots&\ddots&\vdots\\1&1&\dots&1\end{smallmatrix}\right)\text{ and }T^{(1)}:=\frac2{d(d+1)}\left(\begin{smallmatrix}1&1&\dots&1\\0&1&\dots&1\\\vdots&\vdots&\ddots&\vdots\\0&0&\dots&1\end{smallmatrix}\right).$$We denote by $t^{(0)}_{i,j}$ and $t^{(1)}_{i,j}$ their entries respectively, and by $a^{(n)}_{i,j}$ the ones of $\frac{A_n}{\Vert A_n\Vert}$. Note that $\mathcal E$ is a subset of $\cup_k\cap_{\ell\ge k}\mathcal E_{k,\ell}$, where$$\textstyle\mathcal E_{k,\ell}:=\{(A_n)_{n\in\mathbb N}\;;\;\exists i,j\ \exists h\in\{0,1\}\ \vert a^{(2\ell-h)}_{i,j}-t^{(h)}_{i,j}\vert\ge\frac1k\}$$because, if $(A_n)_{n\in\mathbb N}$ is in the complementary of $\cup_k\cap_{\ell\ge k}\mathcal E_{k,\ell}$, the sequence $\frac{A_n}{\Vert A_n\Vert}$ has limit-points $T^{(0)}$ and $T^{(1)}$ without common left-eigenvector.

Given $\ell\in\mathbb N$, there exists $\beta>0$ for the set $\{(A_n)_{n\in\mathbb N}\;;\;\vert a^{(2\ell-h)}_{i,j}-t^{(h)}_{i,j}\vert\ge\frac1k\}$ to have probability $\le1-\beta$ for any $h\in\{0,1\},i,j\in\{1,\dots,d\}$, and $\beta$ does not depend on $\ell$. Hence $\mathcal E_{k,\ell}$ has probability $\le1-\beta^{2d^2}$ and consequently $p^*(\cap_{\ell\ge k}\mathcal E_{k,\ell})=0$, implying $p^*(\mathcal E)=0$.\end{proof}

\begin{corollary}\label{nec}Let $(A_n)_{n\in\mathbb N}$ be a non-eventually constant sequence of matrices that belong to a finite set $\mathcal M\subset\mathcal M_d(\mathbb C)$, and suppose that no couple of elements of $\mathcal M$ has a common left-eigenvector. If $\forall n\ P_n\ne0$, the sequence $\frac{P_n}{\Vert P_n\Vert}$ diverges.
\end{corollary}

\begin{proof}[\bf Proof] If the sequence $A_n$ is non-eventually constant, then $\frac{A_n}{\Vert A_n\Vert}$ is no more eventually constant otherwise $\mathcal M$ contains two collinear matrices whose left-eigenvectors are the same. Hence the sequence $\frac{A_n}{\Vert A_n\Vert}_n$ admits at least two distinct limit-points $\frac{M}{\Vert M\Vert}$ and $\frac{M'}{\Vert M'\Vert}$ with $M,M'\in\mathcal M$ (without common left-eigenvector) and, using Theorem \ref{aediv}, the sequence $\frac{P_n}{\Vert P_n\Vert}$ diverges.\end{proof}

\begin{remark}\label{first}If $A_n$ belongs to a finite set of nonnegative matrices with common left-eigenvector $^tu$ for the eigenvalue 1, one deduces from \cite[Theorem~2.1]{MP} the necessary and sufficient condition for the sequence $P_n$ (or, equivalently, for the sequence $\frac{P_n}{\Vert P_n\Vert}$) to converge to a matrix whose rows are collinear to~$^tu$.\end{remark}

\begin{remark}\label{second}In the example of Subsection \ref{appli}, if the matrices $A_n$ belong to the set $\{M_0,M_1,M_2\}$ defined in Theorem \ref{betareal} and if the sets $\{n\;;\;A_n=M_0\}$ and $\{n\;;\;A_n\in\{M_1,M_2\}\}$ are infinite, the sequence $\frac{P_n}{\Vert P_n\Vert}$ diverges because

$*$ $M_1$ has left-eigenvectors $(0,0,0,0,0,a,b)$, $(1.22,0,1.49,1.82,1,0,0)$ and $(-0.72,0,0.52,-0.38,1,0,0)$ approximately,

$*$ $M_2$ has left-eigenvectors $(a,b,-a,0,-a,c,d)$ and $(a,0,0,-a,b,0,a)$,

$*$ no of them is eigenvector of $M_0$.

%On the other hand the sequence $\frac{P_nc}{\Vert P_nc\Vert}$ converges by Lemma \ref{thnec}.
\end{remark}

\begin{remark}\label{third}At the end of \cite{EF97}, the authors note that Theorem~1 can be seen as a necessary and sufficient condition, that is, given a finite set $\mathcal M=\{A_1,\dots,A_m\}$ of complex-valued matrices without common left-eigenvector for the eigenvalue~1, an infinite product containing an infinite number of each $A_i$ converges iff there exist a (sub-multiplicative) norm $\Vert\cdot\Vert$ and $q<1$ such that\begin{equation}\label{infinite}\textstyle\forall i\ \Vert A_i\Vert\le1\text{ and }\Vert A_{i_1}\cdots A_{i_k}\Vert\le q\text{ if }\forall i\ \exists j\ i=i_j.\end{equation}The authors regret that (\ref{infinite}) requires an infinite number of conditions, if $\mathcal M$ has more than two elements. Nevertheless it seems that \cite{EF97} can be extended in certain cases, for instance if $\mathcal M=\{A_1,A_2,A_3\}$ is a set of three nonnegative matrices: one can reduce (\ref{infinite}) to a finite number of conditions by using the finite set $\mathcal L_i$ of the limit-points of the sequences $n\mapsto{A_i}^n$ (non-empty if $\Vert A_i\Vert\le1$). One obtains that any infinite product containing an infinite number of each $A_i$ converges iff there exist a sub-multiplicative norm $\Vert\cdot\Vert$, a real $q<1$ and an integer $n$ such that:

$*$ $\forall i\ \Vert A_i\Vert\le1$ and $\exists L\in\mathcal L_i\ \Vert{A_i}^n-L\Vert\le\frac{1-q}2$

$*$ $\forall i,j,k$ permutation of $1,2,3$ and $B_j\in\{A_j,\dots,{A_j}^n\}\cup\mathcal L$, $\Vert A_iB_jA_k\Vert\le q$.\end{remark}

\appendix\section{Multifractal formalism \cite{G,HP,Ols,Pesin} and weak-Gibbs property \cite{Y,IY,FO}}\label{mf}

Both notions are related, and the weak-Gibbs property of a sofic measure is related to the uniform convergence of a sequence of the form $\frac{P_nc}{\Vert P_nc\Vert}$. 

\begin{definition}\label{2definitions}(i) Let $\nu$ be a probability measure on the product-set $\Omega_b$. One says that $\nu$ has the {\sl weak-Gibbs property} if there exists a continuous map $\Phi:\Omega_b\to\mathbb R$, called a potential of $\mu$, such that
\begin{equation}\label{defGibbsOmega}
\lim_{n\to\infty}\Big(\frac{\nu([\omega_1\dots\omega_n])}{\exp\big(\sum_{k=0}^{n-1}\Phi(\sigma^k\omega\big)}\Big)^{\frac1n}=1\text{ uniformly on }\Omega_b.
\end{equation}

(ii) Let $\nu$ have the weak-Gibbs property on $\Omega_b$ and let $\mathcal S$ be a system of affine contractions $S_k:[0,1]\mapsto[0,1]$ for
$k\in\{0,\dots,b-1\}$, with $[0,1)=\bigcup_{k\in\{0,\dots,b-1\}}S_k([0,1))$ (disjoint union).
We define the probability measure $\mu$ on the interval $[0,1]$, by
$$\textstyle\mu([\omega_1\dots\omega_n]_\mathcal S):=\nu([\omega_1\dots\omega_n])\text{ for any }[\omega_1\dots\omega_n]_{\mathcal S}:=S_{\omega_1}\circ \cdots\circ S_{\omega_n}([0,1));.$$
One says that $\mu$ has the {\sl weak-Gibbs property with respect to $\mathcal S$, or with respect to the inverse map defined by $T(x)={S_k}^{-1}(x)$ for any $x\in S_k([0,1))$ and $k\in\{0,\dots,b-1\}$}.\end{definition}

The weak-Gibbs property of a probability measure on $[0,1]$ depends on a system of affine contractions, but implies the multifractal formalism, that depends only on the measure. 

The {\it multifractal analysis} is a particular way of analysing the local structure of measures. Let $\mu$ be a probability measure on $[0,1]$; its {\it singularity spectrum} $\tau_{\rm sing}:\mathbb R\to[-\infty,1]$ is defined by
$$
\begin{array}{lcl}\tau_{\rm sing}(\alpha)&:=&\text{H-dim}(E_\mu(\alpha))\quad(\text{Hausdorff dimension of }E_\mu(\alpha)),\text{ where}\\
E_\mu(\alpha)&:=&\{x\;;\;\dim_\mu(x)\text{ exists and }\dim_\mu(x)=\alpha\}\quad(\text{level-set of }\mu),\\
\dim_\mu(x)&:=&\lim_{r\to0}\log(\mu([x-r,x+r]))/\log r\quad(\text{local dimension of }\mu)\end{array}
$$
(by convention, the Hausdorff dimension of the empty set is $-\infty$). Its {\it scale spectrum} (or {\it$L^q$-spectrum}) $\tau_{\rm scale}:\mathbb R\to[-\infty,1]$ is defined by
$$
\tau_{\rm scale}(q):=\liminf_{r\to0}\Big(\log_r\Big(\inf_\mathcal I\Big(\sum_{k=1}^n(\mu(I_k))^q\Big)\Big)\Big),
$$
where $\mathcal I$ is the set of the covers of the support of $\mu$ by closed intervals of length $r$.

The scale spectrum is easier to compute or to approximate than the singularity spectrum. The multifractal formalism, defined as follows, enables to compute the second when one knows the first.\begin{definition}One says that $\mu$ satisfies the {\it multifractal formalism} if the singularity spectrum is the Legendre-transform conjugate of the scale spectrum (see for instance \cite{HJKPS}) in the sense that, for any $\alpha\in\mathbb R$,$$
\tau_{\rm sing}(\alpha)=\inf_{q\in\mathbb R}(\alpha q-\tau_{\rm scale}(q)).
$$\end{definition}

According to \cite[Theorem A']{FO} the weak-Gibbs measures on $[0,1]$ satisfy the multifractal formalism.

Below is a method to prove the weak-Gibbs property for a sofic measure:

\begin{theorem}\label{weak-Gibbs}A sofic measure $\nu$ on $\Omega_b$ has the weak-Gibbs property if the vectors and matrices of its linear representation satisfy the conditions:

$*$ $c_{\omega,n}=\frac{M_{\omega_1}\cdots M_{\omega_n}c}{\Vert M_{\omega_1}\cdots M_{\omega_n}c\Vert}$ converges uniformly on $\Omega_b$ to a limit~$c_\omega$;

$*$ $\forall i\in\{0,\dots,b-1\},\ (r_ic_{\omega,n})^{\frac1n}$ converges to $1$ uniformly on $\Omega_b$;

$*$ $\forall\omega\in\Omega_b\ M_ic_\omega\ne0$.\end{theorem}

\begin{proof}[\bf Proof]We compute
\begin{equation}\label{muandPhi}\begin{array}{rcl}
\frac{\nu([\omega_1\dots\omega_n])}{\Vert c\Vert}&=&\frac{r_{\omega_1}M_{\omega_2}\cdots M_{\omega_n}c}{\Vert M_{\omega_2}\cdots M_{\omega_n}c\Vert}\cdot\frac{\Vert M_{\omega_2}\cdots M_{\omega_n}c\Vert}{\Vert M_{\omega_3}\cdots M_{\omega_n}c\Vert}\cdots\frac{\Vert M_{\omega_n}c\Vert}{\Vert c\Vert}\\&=&r_{\omega_1}c_{\sigma\omega,n-1}\cdot\Vert M_{\omega_2}c_{\sigma^2\omega,n-2}\Vert\cdots\Vert M_{\omega_n}c_{\sigma^n\omega,0}\Vert\end{array}
\end{equation}
and
\begin{equation}\label{Phiandmu}
\exp\big(\sum_{k=0}^{n-1}\Phi(\sigma^k\omega)\big)=\Vert M_{\omega_1}c_{\sigma\omega}\Vert\cdots\Vert M_{\omega_n}c_{\sigma^n\omega}\Vert.
\end{equation}
By the second condition, $(r_{\omega_1}c_{\sigma\omega,n-1})^{\frac1n}$ converges uniformly to $1$, as well as $\Vert M_{\omega_1}c_{\sigma\omega}\Vert^{\frac1n}$ because the third condition and the continuity of $\omega\mapsto c_\omega$ imply $\forall\omega\ \frac1K\le\Vert M_{\omega_1}c_{\sigma\omega}\Vert\le K$ with $K>0$ constant. Since $c_{\omega,n}$ converges uniformly to $c_\omega$, there exists $\varepsilon_n\to0$ such that, for any $\omega\in\Omega_b$,
$$
1-\varepsilon_n\le\frac{\Vert M_{\omega_1}c_{\sigma\omega,n-1}\Vert}{\Vert M_{\omega_1}c_{\sigma\omega}\Vert}\le1+\varepsilon_n\ .
$$
So (\ref{muandPhi}) and (\ref{Phiandmu}) imply
$$\begin{array}{rcl}\frac{r_{\omega_1}c_{\sigma\omega,n-1}}{\Vert M_{\omega_1}c_{\sigma\omega}\Vert}(1-\varepsilon_{n-1})\cdots(1-\varepsilon_1)&\le&\frac{\nu([\omega_1\dots\omega_n])/\Vert c\Vert}{\exp\big(\sum_{k=0}^{n-1}\Phi(\sigma^k\omega)\big)}\\&\le&\frac{r_{\omega_1}c_{\sigma\omega,n-1}}{\Vert M_{\omega_1}c_{\sigma\omega}\Vert}\cdot(1+\varepsilon_{n-1})\cdots(1+\varepsilon_1).\end{array}$$
Since $\lim_{n\to\infty}\big(\prod_{k=1}^{n-1}(1\pm\varepsilon_k)\big)^{\frac1n}=1$, $\nu$ satisfies (\ref{defGibbsOmega}).\end{proof}

\section{Proof of Theorem \ref{cvandGibbs}}\label{proofcvandGibbs}

We consider the matrices $M_0,\dots,M_{k-1}$ defined in Proposition \ref{Fcltsde}, and the integer $1\le r\le k-1$ such that $N=(k-1)q+r$.

%For any word $w=\omega_1\dots\omega_n\in\{0,\dots,k-1\}^n$ we denote by $M_w$ the product-matrix $M_{\omega_1}\cdots M_{\omega_n}$. 
\subsection{Properties of the matrices $M_i$}

\begin{lemma}\label{positive}(i) For any word $\omega_1\dots\omega_q\in\{0,\dots,k-1\}^q$, if $0<\omega_1<r$ the product-matrix $M_{\omega_1}\cdots M_{\omega_q}$ is positive.

(ii) For any word $\omega_1\dots\omega_q\in\{0,\dots,k-1\}^q$, if $\omega_1>r$ the product-matrix $M_{\omega_1}\cdots M_{\omega_q}$ is positive except the last row that is null.

(iii)  For any word $\omega_1\dots\omega_{q+1}\in\{0,\dots,k-1\}^{q+1}$, if $(\omega_1,\omega_2)=(0,r)$ the product-matrix $M_{\omega_1}\cdots M_{\omega_{q+1}}$ is positive.\end{lemma}

\begin{proof}[\bf Proof](i) We use the three $(0,1)$-matrices of same size as the $M_i$, whose nonnull entries are the ones of the the diagonal and the ones just at the left or the right of the diagonal, except in the last row of $Y$ and $Z$ and the first row of $Z$:$$X:=\left(\begin{smallmatrix}1&1&0&\dots&0&0&0\\1&1&1&\dots&0&0&0\\\vdots&\vdots&\vdots&\ddots\vdots&\vdots&\vdots\\0&0&0&\dots&1&1&1\\0&0&0&\dots&0&1&1\end{smallmatrix}\right),\ Y:=\left(\begin{smallmatrix}1&1&0&\dots&0&0&0\\1&1&1&\dots&0&0&0\\\vdots&\vdots&\vdots&\ddots\vdots&\vdots&\vdots\\0&0&0&\dots&1&1&1\\0&0&0&\dots&0&0&0\end{smallmatrix}\right),\ Z:=\left(\begin{smallmatrix}0&0&0&\dots&0&0&0\\1&1&1&\dots&0&0&0\\\vdots&\vdots&\vdots&\ddots\vdots&\vdots&\vdots\\0&0&0&\dots&1&1&1\\0&0&0&\dots&0&0&0\end{smallmatrix}\right).$$

(i) For any $0<i<r$ one has $M_i\ge mX$ with $m=\min_jp_j$. For the other values of $i$, $M_i\ge mZ$. Hence if $0<\omega_1<r$, $M_{\omega_1}\cdots M_{\omega_q}\ge m^qXZ^{q-1}>0$.

(ii) If $\omega_1>r$, $M_{\omega_1}\cdots M_{\omega_q}\ge m^qYZ^{q-1}$. The entries of the $q$ first rows of $YZ^{q-1}$ are positive, and the entries of the last row of $M_{\omega_1}\cdots M_{\omega_q}$ are null.

(iii) $M_0M_r\ge m^2X$ hence  $M_{\omega_1}\cdots M_{\omega_{q+1}}\ge m^{q+1}XZ^{q-1}>0$.\end{proof}%(ii) The ratio $\frac{\sum_jm_{i,j}c_j}{\sum_jm_{i',j}c_j}$, $m_{i,j}$ entries of $M_w$ and $c_j$ entries of $c$, cannot exceed $\Lambda(M_w)$ if $\{j\;;\;m_{i,j}\ne0\}$ is a subset of $\{j\;;\;m_{i',j}\ne0\}$.

%(ii) This follows from the formulas ${M_0}^n=\left(\begin{smallmatrix}{p_0}^n&0\\{p_0}^n\sum_1^n{p_0}^{-i}{M'_0}^{i-1}v_r&{M'_0}^n\end{smallmatrix}\right)$ and ${M'_0}^n=\rho(M'_0)^n(c'_0\ell'_0+N_n)$ with $\lim_{n\to\infty}N_n=0$, $c'_0,\ell'_0$ right and left eigenvectors of $M'_0$ respectively, and the analogous formulas for ${M_r}^n$ and ${M'_r}^n$.

\subsection{Proof of Theorem \ref{cvandGibbs}}Since this theorem requires an uniform convergence, we first give a method to relate pointwise and uniform convergence for the sequences of the form $\frac{P_nc}{\Vert P_nc\Vert}$.

\begin{lemma}\label{pointwisuniform}Given a set of matrices $\mathcal M=\{M_0,\dots,M_{a-1}\}$ and a column-vector $c$ such that the product-matrices $A_1\cdots A_nc$ are nonnull when $A_i\in\mathcal M$, the sequence $ c_{\omega,n}=\frac{M_{\omega_1}\cdots M_{\omega_n}c}{\Vert M_{\omega_1}\cdots M_{\omega_n}c\Vert}$ converges uniformly over $\Omega_a=\{0,\dots,a-1\}^\mathbb N$ if and only if the following pointwise convergence holds for any $\omega\in\Omega_a$:\begin{equation}\label{modifunif}\lim_{n\to\infty}\big(\sup\{\Vert c_{\xi,p}-c_{\xi,p'}\Vert\;;\;\xi_1=\omega_1,\dots,\xi_n=\omega_n,,p,p'\ge n\}\big)=0.\end{equation}

%(ii) In particular if $c$ and the matrices $M_i$ are nonnegative, the sequence $c_{\omega,n}$ converges uniformly over $\Omega_a$ if the normalized columns of $M_{\omega_1}\cdots M_{\omega_n}$ converge pointwise to the same limit.
\end{lemma}

\begin{proof}[\bf Proof]The direct implication is given by the Cauchy criterion. The converse implication follows from a cover of $\Omega_a$ by some suitable cylinder sets $[\omega_1\dots\omega_n]:=\{\xi\in\Omega_a\;;\;\forall k\ \xi_k=\omega_k\}$. Suppose that (\ref{modifunif})~holds for any $\omega\in\Omega_a$. Given $\varepsilon>0$ and $\omega\in\Omega_a$ there exists a rank $n(\omega)$ such that$$\xi\in[\omega_1\dots\omega_{n(\omega)}],p,q\ge n(\omega)\ \Rightarrow\ \left\Vert c_{\xi,p}-c_{\xi,p'}\right\Vert\le\varepsilon.$$Each cylinder set $[\omega_1\dots\omega_{n(\omega)}]$ is an open set for the usual topology and contains $\omega$. Because $\Omega_a$ is compact, it is covered by finitely many of such sets, say $[\omega^{i}_1\dots\omega^{i}_{n(\omega^i)}]$ for $i=1,\dots,N$. Let $\displaystyle n=\max_{i\in\{1,\dots,N\}}n(\omega^i)$. The inequality $\left\Vert c_{\xi,p}-c_{\xi,p'}\right\Vert\le\varepsilon$ holds for any $\omega\in\Omega_a$, $\xi\in[\omega_1\dots\omega_n]$, $p,p'\ge n$, proving that the sequence $c_{\omega,n}$ is uniformly Cauchy and converges uniformly.%(ii) Because $c_{\xi,p}$ and $c_{\xi,p'}$ are in the convex hull of the normalized columns of $M_{\omega_1}\cdots M_{\omega_n}$. %The norm of their difference cannot exceed the largest norm of the difference between two normalized columns of $M_{\omega_1}\cdots M_{\omega_n}$.
\end{proof}

%To prove that (\ref{modifunif}) holds, we first note that the vector $c_{\xi,p}$ is a linear combination of the columns of $c_{\omega,n}$:$$\textstyle c_{\xi,p}=\sum_j\lambda_jc_{\xi,n}u_j\text{ with }\lambda_j=\frac{}{}.$$

%and :=\frac{M_{\omega_1}\cdots M_{\omega_n}c}{\Vert M_{\omega_1}\cdots M_{\omega_n}c\Vert}$

$c_{\xi,p}$ and $c_{\xi,p'}$ are linear nonnegative combinations of the vectors $\frac{M_{\omega_1}\cdots M_{\omega_n}u_i}{\Vert M_{\omega_1}\cdots M_{\omega_n}u_i\Vert}$. When two vectors $x$ and $y$ of a normed space are linear combinations $\sum_i\alpha_iz_i$ of some vectors $z_i$, with $\alpha_i\ge0$ and $\sum_i\alpha_i=1$, one has obviously \begin{equation}\label{maxdist}\textstyle\Vert x-y\Vert\le\max_{i,j}\Vert z_i-z_j\Vert.\end{equation}Let $v$ and $w$ be two normalized nonnegative vectors of $\mathbb R^d$ and suppose that the set $I$ of the indices of their nonnul entries is the same. The projective distance $\delta(v,w):=\max_{i,j\in I}\big(\log\frac{v_iw_j}{v_jw_i}\big)$ \cite[Section 3]{Sen81} is at least $\max_i\big(\log\frac{v_i}{w_i}\big)=\max_i\vert\log v_i-\log w_i\vert$ because there exists necessarily $j$ such that $v_j\le w_j$, and consequently\begin{equation}\label{deltad}\textstyle\delta(v,w)\ge\max_i\vert v_i-w_i\vert\ge\frac1d\Vert v-w\Vert.\end{equation}We also use the Birkhoff contraction coefficient $\tau_\mathcal B$ and both expressions of $\tau_\mathcal B(M)$ at the beginning of \cite[Subsection 3.4]{Sen81}. We deduce, for any nonnegative $d\times d$ matrices $M$ and $M'$,\begin{equation}\label{dist}\textstyle\tau_\mathcal B(M)=\max_{i,j}\tanh\big(\frac{\delta(Mu_i,Mu_j)}4\big)\text{ and }\tau_\mathcal B(MM')\le\tau_\mathcal B(M)\tau_\mathcal B(M').\end{equation}The ineqaulity $\tau_\mathcal B(M)<1$, true for the positive matrices, can be extended to the matrices whose rows are either positive or null, according to the end of the proof of \cite[Theorem 3.10]{Sen81}.

Let the sequence $\omega\in\{0,\dots,k-1\}^\mathbb N$ be neither eventually 0 nor eventually~$r$. Then there exist infinitely many occurrences of the words involved in Lemma~\ref{positive}, that is, the words of length $q$ with first letter distinct from $0$ and $r$, and the words of length $q+1$ whose two first letters are 0 and $r$. Let $\mathcal W$ be this set of words, $\displaystyle\tau:=\max_{h\in\{q,q+1\}\atop\omega_1\dots\omega_h\in\mathcal W}\tau_\mathcal B(M_{\omega_1}\cdots M_{\omega_h})<1$, and let $k(n)$ be the largest number of disjoint intervals $\{i+1,\dots,i+h\}\subset\{1,\dots,n\}$ such that $\omega_{i+1}\dots\omega_{i+h}\in\mathcal W$. The convergence in (\ref{modifunif}) follows because (\ref{maxdist}), (\ref{deltad}) and (\ref{dist}) imply $\Vert c_{\xi,p}-c_{\xi,p'}\Vert\le\tau^{k(n)}\cdot d$.

If the sequence $\omega$ is eventually 0, the Birkhoff contraction coefficient of the powers of $M_0=\left(\begin{smallmatrix}p_0&0\\v_0&M'_0\end{smallmatrix}\right)$ is 1, so the previous method does not apply and a detailed computation is necessary to bound $\Vert c_{\xi,p}-c_{\xi,p'}\Vert$. 

Let $n\in\mathbb N$, with $n\ge n_0$ if there exists $n_0:=$ the largest rank $n$ such that $\omega_n\ne0$. Given $\xi\in[\omega_1\dots\omega_n]$ and $p\ge n$ we use the decomposition$$\textstyle M_{\xi_1}\cdots M_{\xi_p}=PM_0^hQ$$where $P=\prod_{i=1}^{n_0}M_{\omega_i}$ and $Q$ is either the identity matrix or the product of the $M_{\xi_i}$ between the smallest rank $i>n$ such that $\xi_i\ne0$ (if it exixts) and the rank $p$.

Recall that $M_0=\left(\begin{smallmatrix}p_0&0\\v_0&M'_0\end{smallmatrix}\right)$, where $v_0=\left(\begin{smallmatrix}p_k\\\vdots\\p_{qk}\end{smallmatrix}\right)\ne0$. Let $\rho$ (resp. $\rho'$) be the spectral radius of $M_0$ (resp. $M'_0$). There exist $\ell$ (resp. $v$) left-eigenvector (resp. right-eigenvector) of $M_0$ such that$$\textstyle{M_0}^h=h^\alpha\rho^h(v\ell+N_h)\text{ with }\alpha\in\{0,1\}\text{ and }\lim_{h\to\infty}N_h=0.$$We compute $\ell=\left(\begin{smallmatrix}1&0\end{smallmatrix}\right)$ if $p_0\ge\rho'$, $\ell=\left(\begin{smallmatrix}\frac{\ell'v_0}{\rho'-p_0}&\ell'\end{smallmatrix}\right)$ if $p_0<\rho'$, where $\ell'$ is a positive eigenvector of $M'_0$, primitive matrix by Lemma \ref{positive}.

Applying the classical inequality $\big\Vert\frac x{\Vert x\Vert}-\frac y{\Vert y\Vert}\big\Vert\le\frac{2\Vert x-y\Vert}{\Vert y\Vert}$ to the vectors $c_{\xi,p}=\frac{PM_0^hRSc}{\Vert PM_0^hRSc\Vert}$ and $c'_{\xi,p}:=\frac{Ph^\alpha\rho^hv\ell Qc}{\Vert Ph^\alpha\rho^hv\ell Qc\Vert}=\frac{Pv}{\Vert Pv}$ we obtain\begin{equation}\label{cc'}\textstyle\Vert c_{\xi,p}-c'_{\xi,p}\Vert\le\frac{2\Vert PN_hQc\Vert}{\Vert Pv\ell Qc\Vert}.\end{equation}We distinguish two cases, either $Q$ is the product of less than $q$ matrices $M_{\xi_i}$, either $Q=RS$ where $R$ is the product of $q$ matrices $M_{\xi_i}$. In the first case, $Q$ belongs to the finite set of the products of less that $q$ matrices of $\{M_0,\dots,M_{k-1}\}$, hence $\Vert Qc\Vert$ has an upper bound and the real $\ell Qc$ has a positive lower bound, so the convergence in (\ref{modifunif}) follows from (\ref{cc'}) with $\lim_{h\to\infty}N_h=0$.

In the second case, there exist from Lemma~\ref{positive} two reals $0<m\le K$ such that$$\textstyle mw(^tu)\le R\le Ku(^tu)\text{ with }w=\sum_{i=1}^qu_i\text{ and }u=\sum_{i=1}^{q+1}u_i.$$We use the upper bound of $R$ at the numerator of (\ref{cc'}), with $Q=RS$, and the lower bound at the denominator. The convergence in (\ref{modifunif}) follows after simplification by $(^tu)Sc$.

The proof is the same if the sequence $\omega$ is eventually $r$. Note that this fails if $p_N$ is larger than the spectral radius of $M'_r$ and $\frac N{k-1}$ is not an integer, because in this case $M_{r+1}$ belongs to the set of matrices $\{M_0,\dots,M_{k-1}\}$, and the product of the left-eigenvector of $M_r$ (that is, $\left(\begin{smallmatrix}0&\dots&0&1\end{smallmatrix}\right)$) by $M_{r+1}$ is null.

%${M_0}^n=\left(\begin{smallmatrix}{p_0}^n&0\\\frac{{p_0}^n-\rho^n}{p_0-\rho}c\ell&{M_0}^n\end{smallmatrix}\right)+o(\max(p_0}^n,\rho^n))$, where we replace $\frac{{p_0}^n-\rho^n}{p_0-\rho}$ by $n{p_0}^{n-1}$ if $p_0=\rho$. ${M_0}^n=\left(\begin{smallmatrix}{p_0}^n&0\\{p_0}^n\sum_1^n{p_0}^{-i}{M'_0}^{i-1}v_0&{M'_0}^n\end{smallmatrix}\right)$ ${M'_0}^n=\rho(M'_0)^n(c'_0\ell'_0+N_n)$ with $\lim_{n\to\infty}N_n=0$, $c'_0,\ell'_0$ right and left eigenvectors of $M'_0$ respectively, and the analogous formulas for ${M_r}^n$ and ${M'_r}^n$.
%The last row in the matrices of the set $\mathcal M':=\{M_{i_1}\cdots M_{i_q}\;;\;i_1>r\}$ is null. So we can replace, in the proof, the column-vector $u=\sum_1^{q+1}u_i$ by $u'=\sum_1^qu_i$, the row-vector $^tu$ being unchanged. To prove (\ref{modifunif}) we use the pointwise convergence of the sequence $\frac{P_nv}{\Vert P_nv\Vert}$ for $v=\frac{M'u_i}{\Vert M'u_i\Vert}$, M'\in\mathcal M'$.

(ii) Denoting by $[\omega_1\dots\omega_n]$ the set of the reals $x\in[0,1)$ whose $k$-expansion begins by $\omega_1,\dots,\omega_n$, one has by (\ref{linearep})$$\nu_{k,p}([\omega_1\dots\omega_n]+h)={^tu_{h+1}}M_{\omega_1}\cdots M_{\omega_n}c_\nu.$$If $(H)$ holds, Theorem \ref{weak-Gibbs} applies and the measures $\nu_{k,p}(\cdot+h)$ have the weak-Gibbs property on the interval $[0,1]$.

Conversely suppose that the measures $\nu_{k,p}(\cdot+h)$ have the weak-Gibbs property on $[0,1]$. The definition of "$\nu$ weak-Gibbs" implies\begin{equation}\label{omega2n}\lim_{n\to\infty}\left(\frac{\nu([\omega_1\dots\omega_{2n}])}{\nu([\omega_1\dots\omega_n])\nu([\omega_{n+1}\dots\omega_{2n}])}\right)^{\frac1n}=1\ \text{uniformly on }\omega.\end{equation}Consequently $(H)$ holds, otherwise one has $p_N>\rho(M'_r)$ and (\ref{omega2n}) is false for $\nu=\nu_{k,p}(\cdot+h)$ and $\omega_1\dots\omega_{2n}=r^n(r+1)^n$.

\section{Proof of Theorem \ref{betareal}}\label{proofbetareal}We denote by $\mathcal I(X)$ the set of the indices of the non-null entries in the vector or matrix $X$, and we denote by $\mathcal Z(X)$ the $(0,1)$-vector or matrix such that $\mathcal I(\mathcal Z(M))=\mathcal I(M)$. The set of matrices $\mathcal M=\{M_0,M_1,M_2\}$ we consider is defined in Subsection \ref{appli}.

To prove Theorem~\ref{betareal} we do not use the Birkhoff contraction method, because there exists a non-countable set of sequences $(\omega_n)_{n\in\mathbb N}\in\{0,1,2\}^\mathbb N$, like the sequences of the form$$\omega_1\omega_2\dots=(1100000)^{i_1}(1012)^{j_1}(1100000)^{i_2}(1012)^{j_2}\dots\text{ with }i_k,j_k\ge1,$$for any product $P=M_{\omega_m}\cdots M_{\omega_n}$, $1\le m\le n$, to 
 have Birkhoff contraction coefficient 1, since $\mathcal I(P)$ does not have the form $I\times J$ (for instance if $P=(M_1)^2(M_0)^5M_1M_0M_1M_2$
 then $\mathcal I(P)=(\{3,4,5\}\times\{1,5,7\})\cup(\{3,4\}\times\{4\})$). %To prove that $\mu_\beta$ has support $[0,1]$ it is sufficient to find a cover of $[0,1)$ by arbitrarily small intervals of positive measure. For this we chose an arbitrarily large integer $n$ and we consider the intervals $i_0+\frac1{\beta}[\omega_1\dots\omega_n]_{\mathcal S}$ and $i_4+\frac1{\beta}[\omega_1\dots\omega_n]_{\mathcal S}$. They have length less than $\frac1{\beta^n}$ by definition of the maps $S_k$~, and positive measure in consequence of (\ref{translatedcylinders}) because the set of the nonnegative columns-vectors with positive first, third and fifth entries is stable by left-multiplication by $M_k$~, $k\in\{0,1,2\}$. One has $[0,1)=\bigcup_kS_k([0,1))$, hence $[0,1)=\bigcup_{\omega_1\dots\omega_n}[\omega_1\dots\omega_n]_{\mathcal S}$ and consequently$$\begin{array}{rcl}[0,1)&=&\big[0,\frac1\beta\big)\bigcup\big[1-\frac1\beta,1\big)\\&=&\displaystyle\Big(\bigcup_{\omega_1\dots\omega_n}\big(i_0+\frac1{\beta}[\omega_1\dots\omega_n]_{\mathcal S})\Big)\bigcup\Big(\bigcup_{\omega_1\dots\omega_n}\big(i_4+\frac1{\beta}[\omega_1\dots\omega_n]_{\mathcal S})\Big).\end{array}$$
 %\begin{theorem}\label{propertiesofmu}The measure $\mu_\beta$ is weak-Gibbs with respect to the system of affine contractions defined in Lemma \ref{SiM} and to the potential defined on $\Omega_3=\{0,1,2\}^\mathbb N$ by $\Phi(\omega)=\log\Vert M_{\omega_1}c_{\sigma\omega}\Vert$~, where the map $\omega\mapsto c_{\omega}$ is defined in Theorem \ref{weak-Gibbs} and $M_0,M_1,M_2,c$ in Lemma \ref{SiM}. Consequently $\mu_\beta$ satisfies the multifractal formalism.\end{theorem}The proof is given in Appendix \ref{example}.
 So, given a sequence of matrices $(M_{\omega_n})_{n\in\mathbb N}$, $\omega_n\in\{0,1,2\}$, we look for an increasing sequence of integers $(s_k)_{k\in\mathbb N}$ to satisfy the conditions of Theorem~\ref{rankone}. 

\subsection{Choice and properties of an increasing sequence of integers}We read the product $P_n=M_{\omega_1}\cdots M_{\omega_n}$ from the right to the left: either there exists $m$ for $M_{\omega_m}\cdots M_{\omega_n}$ to belong to the following set:

$\mathcal E:=\{M_0M_1{M_0}^\alpha,M_1M_1{M_0}^\alpha,M_2M_1{M_0}^\alpha,M_2{M_0}^\alpha,M_0M_0M_1,M_1M_0M_1,$

$M_2M_0M_1,M_0M_1M_1,M_1M_1M_1,M_2M_1M_1,M_2M_1,M_0M_0{M_2}^\alpha,$

$M_0M_0M_1M_0{M_2}^\alpha,M_1M_0M_1M_0{M_2}^\alpha,M_2M_0M_1M_0{M_2}^\alpha,M_1M_1M_0{M_2}^\alpha,$ 

$M_2M_1M_0{M_2}^\alpha,M_2M_0{M_2}^\alpha,M_0M_1{M_2}^\alpha,M_1M_1{M_2}^\alpha,M_2M_1{M_2}^\alpha\}_{\alpha\in\mathbb N}$,

either $M_{\omega_1}\cdots M_{\omega_n}$ is the end of such a product, that is, $M_{\omega_1}\cdots M_{\omega_n}$ belongs to the set

$\mathcal E':=\{I,{M_0}^\alpha,M_1{M_0}^\alpha,M_1,M_0M_1,M_1M_1,{M_2}^\alpha,M_0{M_2}^\alpha,M_1M_0{M_2}^\alpha,$

$M_0M_1M_0{M_2}^\alpha,M_1{M_2}^\alpha\}_{\alpha\in\mathbb N}$.

By reiterating one obtains\begin{equation}\label{matrices e}\textstyle P_n=E_0E_1\cdots E_k\text{ with }k\ge0,\ E_0\in\mathcal E',\ E_1,\dots,E_k\in\mathcal E.\end{equation}

Moreover we prove that for any sequence $(M_{\omega_n})_{n\in\mathbb N}$ not eventually $M_0$ nor eventually $M_2$, there exists $(E_k)_{k\ge0}$ such that\begin{equation}\label{matrices E}\textstyle\forall k\ \exists n_k\ E_0\cdots E_k=P_{n_k}.\end{equation}Indeed

$*$ the matrix $E_0$ involved in (\ref{matrices e} takes the same value for infinitely many $n$ because the exponent $\alpha$ involved in $E_0$ has a bound that depends only on the sequence $(M_{\omega_n})_{n\in\mathbb N}$

$*$ and, by induction, if $E_0,\dots,E_{i-1}$ are already chosen, the exponent $\alpha$ in $E_i$ has a bound that depends only on $(M_{\omega_n})_{n\in\mathbb N}$, and $E_i$ takes infinitely times the same value.

We look at the properties of the matrices of $\mathcal E$, using the following coefficient:%The coefficient $\Lambda(M)$ below is the maximum-ratio of the distances between the column-vectors of the matrix $M$, seen as positive vectors in a space generated by $d'\le d$ canonical vectors, and the hyperplanes generated by $d'-1$ vectors.
 
 \begin{definition}For any square or rectangular matrix $A\ne0$ we denote by $\Lambda(A)$ the maximum-ratio between two entries of same column-index:$$\textstyle\Lambda(A):=\max_{i,i',j\atop a_{i,j}\ne0}\frac{\vert a_{i',j}\vert}{\vert a_{i,j}\vert}.$$\end{definition}Looking for a bound for $\Lambda(AB)$ when $A$ and $B$ are nonnegative, we obtain the following\begin{lemma}\label{lambdaab}Let $\lambda(A):=\max_{i,i'}\frac{\Vert^tu_{i'}A'\Vert}{\Vert^tu_iA\Vert}$, where $A'$ is the submatrix of $A$ whose columns are the ones of $A$, except the columns $Au_j$ such that $\mathcal I(Au_j)=\cup_i\mathcal I(Au_i)$ if they exist. One has \begin{equation}\label{formulaab}\textstyle\Lambda(AB)\le\Lambda(A)+\lambda(A)\Lambda(B).\end{equation}\end{lemma}

\begin{proof}

Let $p_{i,j}$ be the entries of $AB$, and $i,i',k\in\{1,\dots,7\}$ with $p_{i,k}\ne0$. Using the notatio $J_i=\{ j\;;\;a_{i,j}\ne0$ and denoting by $j(k)$ an index such that $b_{j(k),k}=\min_jb_{j,k}$,$$\begin{array}{rcl}\frac{p_{i',k}}{p_{i,k}}&=&\frac{\sum_{j\in J_i}a_{i',j}b_{j,k}+\sum_{j\not\in J_i}a_{i',j}b_{j,k}}{\sum_ja_{i,j}b_{j,k}}\\&\le&\frac{\Lambda(A)\sum_{j\in J_i}a_{i,j}b_{j,k}}{\sum_ja_{i,j}b_{j,k}}+\frac{\Lambda(B)b_{j(k),k}\sum_{j'\not\in J_i}a_{i',j'}}{b_{j(k),k}\sum_ja_{i,j}}.\end{array}$$The entry $a_{i',j'}$, in the numerator of the second term, is an entry of the ${j'}^{\rm th}$ column of $A$ and, on this same column, $a_{i,j'}=0$. It is a column of $A'$ because the hypothesis $p_{i,k}\ne0$ implies $\exists j\ a_{i,j}\ne0$. The inequality (\ref{formulaab}) follows. 

\end{proof}

We explain below why $\Lambda$ is bounded on the set of the products of elements of~$\mathcal E$. The following lemma is obtained by a simple computation, using the formulas ${M_0}^{4\alpha}=\left(\begin{smallmatrix}1&0&0&0&0&0&0\\\alpha&0&1&0&0&0&0\\\alpha&0&0&1&1&0&0\\0&0&0&0&0&0&0\\\alpha&0&0&0&0&0&1\\\alpha&0&0&0&1&0&0\\\alpha&1&0&0&0&0&0\end{smallmatrix}\right)$ and ${M_2}^\alpha=\left(\begin{smallmatrix}1&0&0&0&\alpha&0&1\\0&0&0&0&0&0&0\\1&0&0&0&\alpha-1&0&1\\0&0&0&1&\alpha&0&0\\0&0&0&0&1&0&0\\0&0&0&0&0&0&0\\0&&0&0&0&0&0\end{smallmatrix}\right)$:

\begin{lemma}\label{thelambdas}(i) The maximum of $\Lambda$ (resp. $\lambda$) on $\mathcal E$ is $\Lambda(M_2{M_0}^4)=3$ (resp. $\lambda(M_0M_1M_0)=3$).

(ii) Any product of three matrices of $\mathcal M$ has at least one column larger than $u_1+u_3+u_5$, the first or third or fifth column.%, in particular the one that depends on the exponent $\alpha$.

(iii) The set of the vectors $v\ge u_1+u_3+u_5$, is stable by left-multiplication by $M_0$, $M_1$ or $M_2$.

(iv) For $P\in\{(M_0)^2,M_1M_0,(M_1)^2,M_2\},v\ge u_1+u_3+u_5$ and $u=\sum_{i=1}^7u_i$,$$\textstyle\mathcal Z(Pv)=\mathcal Z(Pu)\in\Big\{\left(\begin{smallmatrix}1\\1\\1\\0\\1\\1\\1\end{smallmatrix}\right),\left(\begin{smallmatrix}1\\1\\1\\1\\1\\0\\0\end{smallmatrix}\right),\left(\begin{smallmatrix}1\\0\\1\\1\\1\\0\\0\end{smallmatrix}\right),\left(\begin{smallmatrix}1\\1\\1\\0\\1\\1\\0\end{smallmatrix}\right)\Big\}.$$Consequently any % $P=M_{\omega_1}M_{\omega_2}M_{\omega_3}$ and any
$P$ product of at least three matrices of $\mathcal M$ has also this property.

(v) For $M\in\mathcal M$ and $v\in\mathcal S:=\{u_1,u_2,u_3,u_4,u_5,u_6,u_7,u_1+u_2+u_5,u_1+u_3,u_1+u_3+u_4,u_1+2u_3+u_4,u_1+u_4+u_5,u_1+u_5,u_2+u_3,u_2+u_3+u_6,u_2+u_7,u_3+u_4,u_3+u_5+u_6,u_3+u_6,u_5+u_7\}$, either $Mv\in\mathcal S$ or $Mv\ge u_1+u_3+u_5$.\end{lemma}

According tothe following graph, whose initial state is the set of vectors $v\ge u_1+u_3+u_5$ and final state the set of vectors $v\ge2(u_1+u_3+u_5)$, any product-matrix $M_{\omega_1}^{i_1}\cdots M_{\omega_k}^{i_k}$ with $k\ge9$ and $\forall j\ \omega_j\ne\omega_{j+1},i_j\ge1$ transforms (by left-multiplication) any vector of the initial state to a vector of the final state, except if $M_1M_0M_0M_1M_0M_0$ is a part of this product:

\includegraphics[scale=0.40]{M0M1M2.jpg}

In view of (\ref{matrices E}), if $(M_{\omega_n})_{n\in\mathbb N}$ is not eventually $M_0$ nor eventually $M_2$, we associate to any $n\ge n_0$ the integer $k$ such that $n_k\le n<n_{k+1}$ and we have\begin{equation}\label{ME}\textstyle P_n=M_{\omega_1}\cdots M_{\omega_n}=E_0\cdots E_kP_{n_k,n}\end{equation}where $P_{n_k,n}$ (product of the $M_{\omega_i}$ for $n_k<i\le n$) either belongs to $\mathcal E$ or is a product of at most 4 matrices of $\mathcal M$.

$*$ From Lemma \ref{thelambdas}(ii) and (iv), for any $n\ge6$ there exists $j=j(\omega,n)\in\{1,3,5\}$ such that $\mathcal Z(P_{n-6,n}u_j)=\mathcal Z(P_{n-6,n}u)$. Consequently $P_{n-6,n}u\le KP_{n-6,n}u_j$, where $K$ is 7 times the maximum-entry of all the products of 6 matrices of $\mathcal M$. Multiplying by $P_{n-6}$ one obtains\begin{equation}\label{Ktimes}P_nu\le KP_nu_j.\end{equation}

$*$ Moreover the above property of the graph implies that all the entries of $P_nu_j$ are at least $\kappa\rho^k$ with $\kappa>0$ and $\rho>1$ constants, except if there exists~$n$ such that $\omega_{n+1}\dots\omega_{n+6}=100100$.

$*$ From the last item of Lemma \ref{thelambdas}, for any $i$ either $P_{3,n}u_i\ge u_1+u_3+u_5$ or $P_{3,n}u_i\in\mathcal S$, hence either $\mathcal Z(P_nu_i)=\mathcal Z(P_nu)$ or the entries of $P_nu_i$ do not exceed $2\cdot 7^3$.

So the coefficient $\lambda$ defined in Lemma \ref{lambdaab} satisfies the inequality\begin{equation}\label{lambdabound}\textstyle\lambda(P_n)\le C\rho^{-k}\ (C\text{ constant}).\end{equation}Also true if $\exists n\ \omega_{n+1}\dots\omega_{n+6}=100100$, $\lambda(P_n)=0$ because $P_n$ has rank 1.

From Lemma \ref{lambdaab},$$\begin{array}{rcl}\Lambda(P_n)&\le&\Lambda(P_{n_k})+\lambda(P_{n_k})\Lambda(P_{n_k,n})\\&\le&\Lambda(P_{n_{k-1}})+\lambda(P_{n_{k-1}})\Lambda(E_k)+\lambda(P_{n_k})\Lambda(P_{n_k,n})\\&\le&\dots\\&\le&\Lambda(P_0)+\sum_{i=0}^{k-1}\lambda(E_0\cdots E_i)\Lambda(E_{i+1})+\lambda(E_0\dots E_k)\Lambda(P_{n_k,n}).\end{array}$$$\lambda(E_0\cdots E_i)\le C\rho^{-i}$ and $\Lambda(E_{i+1})\le3$, and $P_{n_k,n}$ (if $\not\in\mathcal E$) is the product of at most 4 matrices of $\mathcal M$, so $\Lambda(P_n)$ is bounded and the bound only depends on $E_0=P_{n_0}$ (for instance $\Lambda({M_2}^{n_0})=n_0$).

\subsection{Proof of Theorem \ref{betareal}}We first prove that the sequence $v_{\omega,n}:=\frac{P_nv}{\Vert P_nv\Vert}$ converges pointwise:

$*$ We replace the sequence $(n_k)_{k\in\mathbb N}$ by a subsequence $(s_k)_{k\in\mathbb N}$, for $\frac{P_{s_k}}{\Vert P_{s_k}\Vert}$ to converge in the compact set of the $7\times7$ matrices of norm 1, and for $P_{s_k,s_{k+1}}$ to be the product of at least $k$ matrices of $\mathcal E$. 

$*$ To prove that $(H_1)$ holds, we use the index $j\in\{1,3,5\}$ and the constant $K$ involved in (\ref{Ktimes}), to give a lower bound of$$\textstyle\frac{\Vert P_n\Vert}{{\Vert P_{s_k}\Vert}\Vert P_{s_k,n}\Vert}\ge\frac{\Vert P_{s_k}u_j\Vert\Vert ^tu_jP_{s_k,n}\Vert}{{\Vert P_{s_k}\Vert}\Vert P_{s_k,n}\Vert}\ge K^{-1}\frac{\Vert ^tu_jP_{s_k,n}\Vert}{\Vert P_{s_k,n}\Vert}.$$Since $\Vert P_{s_k,n}\Vert=\sum_i(^tu_i)P_{s_k,n}u\le7\max_i(^tu_i)P_{s_k,n}u$, the r.h.s. is at least $K^{-1}\frac1{\Lambda(P_{s_k,n}u)}$, and $\Lambda(P_{s_k,n}u)\le\Lambda(P_{s_k,n})+\lambda(P_{s_k,n})$ by Lemma \ref{lambdaab}. Recall that $s_{k+1}\le n<s_{k+2}$ and that $P_{s_k,s_{k+1}}$ is a product of matrices of $\mathcal E$. Hence $\Lambda(P_{s_k,n})$ is bounded, as well as $\lambda(P_{s_k,n})$ by (\ref{lambdabound}), and $\frac{\Vert P_n\Vert}{{\Vert P_{s_k}\Vert}\Vert P_{s_k,n}\Vert}$ has a positive lower bound.

$*$ To prove that $(H_2)$ holds, the matrix $P_{s_k,s_{k+1}}$ has some entries $\ge\kappa\rho^k$, so by Lemma \ref{thelambdas}(iv) and (v) the columns $\frac{P_{s_k,s_{k+1}}u_j}{\Vert P_{s_k,s_{k+1}}\Vert}$ such that $\mathcal Z(P_{s_k,s_{k+1}}u_j)\ne\mathcal Z(P_{s_k,s_{k+1}}u)$ tend to 0 when $k\to\infty$ and, if $P$ is a limit-point of $\frac{P_{s_k,s_{k+1}}}{\Vert P_{s_k,s_{k+1}}\Vert}$, the set $\mathcal I(P)$ has the form $I\times J$, implying  $(H_2)$. 

    For any nonnegative vector $v\ge u_1+u_3+u_5$ the sequence $v_{\omega,n}:=\frac{P_nv}{\Vert P_nv\Vert}$ converges to a limit $\ell_\omega$ independant of $v$, because the condition \hbox{$\inf_n\frac{\Vert P_nv\Vert}{\Vert P_n\Vert}>0$} of Theorem \ref{rankone} is satisfied, in consequence of (\ref{Ktimes}). Moreover given $0<\alpha<\frac13$, this convergence is uniform with respect to $v$ if $v$ belongs to the set$$\textstyle\mathcal V_\alpha:=\{v\;;\;v\text{ nonnegative, }\Vert v\Vert=1,\ v\ge\alpha(u_1+u_3+u_5)\},$$ because $V_\alpha$ is a subset of the convex hull of the vectors $e_i=\alpha(u_1+u_3+u_5)+(1-3\alpha)u_i$, according to the relation $v=\sum_i\frac{w_i}{1-3\alpha}e_i$ with $w=v-\alpha(u_1+u_3+u_5)$.

This remains true if  $\forall n>N\ \omega_n=0$ (resp. $\forall n>N\ \omega_n=2$): $v_{\omega,n}$ converges uniformly on $\mathcal V_\alpha$ to $\frac{P_N(u_2+u_3+u_5+u_6+u_7)}{\Vert P_N(u_2+u_3+u_5+u_6+u_7)\Vert}$ (resp. $\frac{P_N(u_1+u_3+u_4)}{\Vert P_N(u_1+u_3+u_4)\Vert}$). So, for any $\omega\in\{0,1,2\}^\mathbb N$,\begin{equation}\label{unifomega}\textstyle\forall\varepsilon>0\ \exists n_\varepsilon\ (v\in\mathcal V_\alpha\text{ and }n\ge n_\varepsilon)\Rightarrow\big\Vert v_{\omega,n}-\ell_\omega\big\Vert\le\varepsilon.\end{equation}We use Lemma \ref{pointwisuniform}: it remains to prove (\ref{modifunif}) or, equivalently, to prove that\begin{equation}\label{unifxi}\textstyle\forall\varepsilon>0\ \exists m_\varepsilon\ (\xi\in[\omega_1\dots\omega_{m_\varepsilon}]\text{ and }p\ge m_\varepsilon)\Rightarrow\Vert c_{\xi,p}-\ell_\omega\Vert\le\varepsilon.\end{equation}

Let $p$ and $m_\varepsilon$ be two integers such that $p\ge m_\varepsilon\ge n_\varepsilon$, and $\xi$ a sequence that coincides with $\omega$ until the rank $m_\varepsilon$. We use the decomposition of (\ref{matrices e}): $M_{\xi_1}\cdots M_{\xi_p}=E_0\cdots E_\ell$. There exist $i,j,h$, depending on $\xi$ and $p$, such that\begin {equation}\label{mnh}i\le n_\varepsilon\le j\text{ and }M_{\xi_i}\cdots M_{\xi_j}=E_h.\end{equation}So if the sequences $\xi$ and $\omega$ coincide until the rank $j$, we obtain (\ref{unifxi}) by applying (\ref{unifomega}) to the vector $v=\frac{E_{h+1}\cdots E_kc}{\Vert E_{h+1}\cdots E_kc\Vert}$ (or $\frac c{\Vert c\Vert}$ if $h=k$), which belongs to $\mathcal V_{\alpha}$ for some value of $\alpha$ because $\Lambda(P)$ is bounded for any $P$ products of elements of $\mathcal E$.

It remains to chose $m_\varepsilon$ for $\xi$ and $\omega$ to coincide necessarily until the rank $j$: if the sequence $\omega$ is not eventually 0 or eventually 2, we chose $m_\varepsilon$ for the word $\omega_{n_\varepsilon}\dots\omega_{m_\varepsilon}$, seen as a concatenation of $0^i$, 1 and $2^i$, to contain at least 6 such words. Since the elements of $\mathcal E$ have at most 5 occurrences, the integer $j$ in (\ref{mnh}) is less than $m_\epsilon$ and the sequences $\xi$ and $\omega$ coincide until the rank~$j$. This is also true if the sequence $\omega$ is eventually 0 (resp. eventually~2): we first chose the integer $n_\varepsilon$ in (\ref{unifomega}) in such a way that $\omega_n=0$ (resp. 2) for any $n\ge n_\varepsilon$, and we chose $m_\varepsilon=n_\varepsilon+2$. So $\omega_{n_\varepsilon}\dots\omega_{m_\varepsilon}=0^3$ and, since ${M_0}^3$ can be only at the end of an element of $\mathcal E$, the integer $j$ in (\ref{mnh}) is at least $m_\epsilon$ and the sequences $\xi$ and $\omega$ coincide until the rank $j$.

(ii) The weak-Gibbs property of $\nu_{\beta,2}\big(\frac\cdot{\beta-1}\big)$ follows, by using Theorem~\ref{weak-Gibbs}.

 %smallest integer $p\ge n_{\varepsilon}$ such that$$\exists h\ M_{\xi_1}\cdots M_{\xi_p}=E_0\cdots E_h.$$From the hypothesis on $n'_\varepsilon$, if $\omega$ is not eventually 0 nor eventually 2, $p$ cannot be larger than $n'_\varepsilon$, so $\xi_i=\omega_i$ for any $i\le p$. If 

%It remains to find $0<\alpha<\frac13$ for $v:=\frac{E_{h+1}\cdots E_kc}{\Vert E_{h+1}\cdots E_kc\Vert}$ to belong to $\mathcal V_\alpha$. The first, third and fifth entry of $v$ are $\ge\frac1{7\Lambda(E_{h+1}\cdots E_kc)}$. By Lemmas \ref{Appl1lL}(ii) and \ref{lL}(i) there exists a constant $\alpha>0$ such that $\frac1{7\Lambda(E_{h+1}\cdots E_kc)}\ge\alpha$, and (\ref{unifomega}) implies~(\ref{unifxi}).\end{proof}

%The proof of the weak-Gibbs property of the measure $\nu\big(\frac \cdot{\beta-1}\big)$ is the same as in the integral case  \S\ref{Mslr}.\end{proof}

%To check the second condition of Theorem \ref{weak-Gibbs} we bound $\Lambda(c_{\omega,n})=\Lambda(P_nc)$. We use Lemma \ref{mathcalI} and apply Lemma \ref{lL}: since $S_0$ is the product of at most three matrices of the set $\{A,B,C\}$ by a power of $A$ or $C$, and by a power of $B$, since $\Lambda(M^n)$ and $\lambda(M^n)$ are in $O(n)$ for $M=A$, $B$ or $C$, $\Lambda(P_nc)=O(n)$ and the second condition holds. The third condition also holds because $\mathcal I(c_\omega)$ contains $\{1,3,5\}$ or, if $\forall n\ \omega_n=0$ (resp. if $\forall n\ \omega_n=2$), $\mathcal I(c_\omega)=\{2,3,5,6,7\}$ (resp. $\{1,3,4\}$).\end{proof}

\section{Triangular matrices}\label{tm}

\subsection{The block-triangular form}\label{tptc}

\begin{definition}Let

$$\begin{array}{lcl}%\mathcal I(A)&:=&\text{ the set of the indices of the nonnull entries in the matrix $A$}
%\\
%\mathcal Z(A)&:=&\text{ the $(0,1)$-matrix such that }\mathcal I(\mathcal Z(A))=\mathcal I(A)
%\\
\mathcal N(A)&:=&\text{ the number of distinct column of $\mathcal Z(A)$}
\\
\mathcal T_\mathcal J&:=&\text{the set of the block-triangular $d\times d$ matrices}\\&&\text{with respect to the partition $\mathcal J=\{J_1,\dots,J_\delta\}$ of $\{1,\dots,d\}$},\\&&(\text{meaning that }\mathcal I(A)\subset\cup_k((J_1\cup\dots\cup J_k)\times J_k))\\
\mathcal T_\mathcal J^1&:=&\text{the set of the matrices }A\in\mathcal T_\mathcal J\text{ such that }\mathcal N(B_{i,j})=1\\&&\text{for any block }B_{i,j}:=(a_{i',j'})_{(i',j')\in J_i\times J_j}.\end{array}$$
\end{definition}

%Given a sequence $\mathcal A=(A_n)_{n\ge1}$ of $d\times d$ matrices we use the notations
%$$\begin{array}{lcl}A_{i,j}&:=&\left\{\begin{array}{ll}A_iA_{i+1}\cdots A_j&(i\le j)\\\text{the identity matrix }I_d&(i>j)\end{array}\right.\\P_n&:=&A_{1,n}=A_1A_2\cdots A_n.\end{array}$$

\begin{theorem}\label{triangular}Let $\mathcal A=(A_n)_{n\ge1}$ be a sequence of nonnegative $d\times d$ matrices. There exist an increasing sequence of nonnegative integers $(n_k)_{k\ge1}$ and a partition $\mathcal J=\mathcal J_\mathcal A$ (unique if we choose $(n_k)_{k\ge1}$ minimal for the lexicographic order) such that, for any $1\le k<\ell$, $P_{n_k,n_\ell}\in T^1_\mathcal J$ and $\mathcal Z(P_{n_k,n_\ell})$ does not depend on $(k,\ell)$.

%(ii) The existence of a vector  $V_\mathcal A:=\lim_{n\to\infty}\frac{P_nv}{\Vert P_nv\Vert}$ independent of the positive column-vector $V$, only depends on $M_\mathcal A$. 
\end{theorem}

\begin{proof}[{\bf Proof}]Let $\delta=\delta_\mathcal A:=\lim_{n\to\infty}\limsup_{n'\to\infty}\mathcal N(P_{n,n'})$.\begin{lemma}\label{delta}There exists $n_1\in\mathbb N$ such that, for any $n\ge n_1$,

$\begin{array}{ll}\exists H_n\text{ infinite subset of }\{n,n+1,\dots\}\ &\mathcal N(P_{n,n'})=\delta\ (n'\in H_n)\\&\mathcal N(P_{n,n'})\le\delta\ (n'\ge\min H_n).\end{array}$
\end{lemma}

\begin{proof}[{\bf Proof}]We note first that $\mathcal N(P_{n,n'})\le\mathcal N(P_{n+1,n'})$ because, if $\mathcal Z(P_{n,n'}u_{j_1})$, $\mathcal Z(P_{n,n'}u_{j_2}),\dots$ are some distinct columns of $\mathcal Z(P_{n,n'})$, then since $P_{n,n'}=A_{n+1}P_{n+1,n'}$ the columns $\mathcal Z(P_{n+1,n'}u_{j_1}),\mathcal Z(P_{n+1,n'}u_{j_2}),\dots$ are necessarily distinct.

Consequently the sequence $n\mapsto\limsup_{n'\to\infty}\mathcal N(P_{n,n'})$ is non-decreasing. It is constant from a rank $n=n_1$ and Lemma \ref{delta} follows.\end{proof}

%\begin{lemma}\label{colAB}If $\mathcal I(AB)=\mathcal I(A)$ and $\mathcal N(B)\le\mathcal N(A)$, then $A$ is upper-block-triangular, with respect to a partition that only depends on the natural partition of $\mathcal I(A)$.\end{lemma}
%\begin{proof}[{\bf Proof}]

%Suppose $\mathcal I(AB)=\mathcal I(A)=\cup_{k=1}^\delta (I_k\times J_k)$, where the $I_k$ are distinct and the $J_k$ make a partition of $\{1,\dots,d\}$.  By reiterating on obtains an indexation such that $k<k'\Rightarrow I_k\not\supset I_{k'}$.

%If $j_1\in J_1,\dots, j_\delta\in J_\delta$, one has $\mathcal I(ABu_{j_1})=I_1,\dots,\mathcal I(ABu_{j_\delta})=I_\delta$ hence $\mathcal I(Bu_{j_1}),\dots,\mathcal I(Bu_{j_\delta})$ are distinct. If there exist two elements $j,j'$ in a set $J_k$ such that $\mathcal I(Bu_j)\ne\mathcal I(Bu_{j'})$, one has $\mathcal N(B)>\delta$ and this contradicts the hypothesis. Hence $\mathcal I(B)=\cup_{k=1}^\delta(I'_k\times J_k)$. The set $I'_k$ cannot intersect any $J_{k'}$ with $k'>k$, otherwise for $i\in I'_k\cap J_{k'}$ and $j\in J_k$ one has $ABu_j\ge b_{i,j}Au_i$ and $\mathcal I(ABu_j)=I_k$ would contain $\mathcal I(Au_i)=I_{k'}$. This proves that $I'_k\subset J_1\cup\dots\cup J_k$ and $B$ is upper-block-triangular with respect to the partition $J_1,\dots,J_\delta$.\end{proof}

We chose, among the elements of $H_{n_1}$, some integers $n_2<n_3<\dots$ such that $\mathcal I(P_{n_1,n_2})=\mathcal I(P_{n_1,n_3})=\dots$. Using the last inequality of Lemma \ref{delta}, we can chose $n_2,n_3,\dots$ large enough to have $\mathcal N(P_{n_k,n_\ell})\le\delta$ when $1\le k<\ell$.

There exist $I_1,\dots,I_\delta$ (distinct), $J_1,\dots,J_\delta$ (partition of $\{1,\dots,d\}$) such that$$\textstyle\forall k>1\ \mathcal I(P_{n_1,n_k})=\cup_{h=1}^\delta (I_h\times J_h).$$At least one of the $I_h$ does not contain any $I_{h'}$ with $h'\ne h$. Since $P_{n_0,n_\ell}=P_{n_0,n_k}P_{n_k,n_\ell}$ and since $\mathcal I(P_{n_0,n_k})=\mathcal I(P_{n_0,n_\ell})$, we deduce that the matrix $P_{n_k,n_\ell}$ cannot have nonnull elements with column-index in $J_h$ and row-index in the complementary of $J_h$. So $\mathcal Z(P_{n_k,n_\ell})$ has the form $P\left(\begin{smallmatrix}B_{k,\ell}&C_{k,\ell}\\0&D_{k,\ell}\end{smallmatrix}\right)P^{-1}$ with $P$ permutation matrix, $\mathcal N(B_{k,\ell})=1$ and $\mathcal N(D_{k,\ell})=\delta-1$ (if $\delta=1$, $B_{k,\ell}$ is the matrix $\mathcal Z(P_{n_k,n_\ell})$ itself).

Assuming that the theorem is true for the sequences of matrices $\mathcal A$ such that $\delta_\mathcal A=\delta-1$, we can suppose $D_{k,\ell}=D$ constant, after replacing $(n_k)_{k\in\mathbb N}$ by a suitable subsequence.

Denoting by $B_k$ (resp. $C_k$) the matrix $B_{k,k+1}$ (resp. $C_{k,k+1}$), one has $\left(\begin{smallmatrix}B_{k,\ell}&C_{k,\ell}\\0&D\end{smallmatrix}\right)=\mathcal Z(\prod_{i=k}^{\ell-1}\left(\begin{smallmatrix}B_i&C_i\\0&D\end{smallmatrix}\right))$ and, since $\mathcal N(B_k)=1$, $B_{k,\ell}=B_k$ except if one of the $B_i$ is null. Replacing $(n_k)_{k\in\mathbb N}$ by a suitable subsequence, $B_{k,\ell}=B_k=0$ in this last case.

$C_{k,\ell}=\mathcal Z(C_kD+B_k\sum_{k<i<\ell-1}C_iD+B_kC_{\ell-1})$. Since the sequence $\ell\mapsto\mathcal Z(\sum_{k<i<\ell-1}C_i)$ is not-decreasing (the entries are not-decreasing), it converges to a $(0,1)$-matrix that we denote by $\mathcal Z(\sum_{k<i<\infty}C_i)$, by abuse of notation. Now $k\mapsto\mathcal Z(\sum_{k<i<\infty}C_i)$ is not-increasing and there exists $\kappa$ such that this sequence converges to $\mathcal Z(\sum_{\kappa<i<\infty}C_i)=:S$.

There exists $K$, infinite set of integers, such that $(B_k,C_k,C_{k-1})=\text{constant}$

$=:(B,C,C')$, so, replacing $(n_k)_{k\ge1}$ by a suitable subsequence one has

$\mathcal Z(P_{n_k,n_\ell})=P\left(\begin{smallmatrix}B&\mathcal Z(CD+BSD+BC')\\0&D\end{smallmatrix}\right)P^{-1}$ and $P_{n_k,n_\ell}\in T^1_\mathcal J$ ($1\le k<\ell$).

To prove the existence of a minimal sequence $(n_k)_{k\ge1}$ for the lexicographical order, that satisfies both conditions of Theorem \ref{triangular}, we consider more generally a set $\mathcal S\ne\emptyset$ of sequences $(n_k)_{k\ge1}$ defined by a condition $\mathcal C$, and we assume the equivalence: $(n_k)_{k\ge1}$ satisfies $\mathcal C$ if and only if the finite sequence $(n_k)_{k=1}^\ell$ satisfies $\mathcal C$ for any $\ell\ge1$. Considering the set $\mathcal S_u$ of the sequences less or equal to a fixed element $u=(n_k)_{k\ge1}\in\mathcal S$, one can easily define by induction the element of $\mathcal S$ minimal for the lexicographical order.\end{proof}

\subsection{Upper triangular $3\times3$ matrices}\label{3}Let $A_n=\left(\begin{smallmatrix}a_n&b_n&c_n\\0&d_n&e_n\\0&0&f_n\end{smallmatrix}\right)$ belong to a finite set of nonnegative upper-triangular matrices with positive diagonal entries, and let $v$ be a positive 3-dimensional column-vector.$$P_n=\left(\begin{matrix}a^*_n&d^*_ns_n&f^*_n(t_n+\tau_n)\\0&d^*_n&f^*_nu_n\\0&0&f^*_n\end{matrix}\right)\text{ with }\left\{\begin{array}{l}{\scriptstyle a^*_n:=a_1\dots a_n,\ d^*_n:=d_1\dots d_n,\ f^*_n:=f_1\dots f_n}\\{\scriptstyle s_n:=\sum_1^n\frac{a^*_{k-1}b_k}{d^*_k},\ t_n:=\sum_1^n\frac{a^*_{k-1}c_k}{f^*_k},}\\{\scriptstyle u_n:=\sum_1^n\frac{d^*_{k-1}e_k}{f^*_k},\ \tau_n:=\sum_2^n\frac{s_{k-1}d^*_{k-1}e_k}{f^*_k}.}\end{array}\right.$$Note that the sequences $n\mapsto s_n,t_n,u_n,\frac{\tau_n}{u_n}$ are non-decreasing (the last sequence because $\frac{\tau_{n+1}}{u_{n+1}}-\frac{\tau_n}{u_n}=\frac{d^*_ne_{n+1}}{f^*_{n+1}u_{n+1}}(s_n-\frac{\tau_n}{u_n})\ge0$). Denoting by $s,t,u,\tau$ their respective limits, $\tau=s$ if $u=\infty$.

$1^{\rm st}$ case: $s$ and ($u$ or $\tau$) are infinite, in this case the limit-points of $\frac{P_n}{\Vert P_n\Vert}$ have the form $\left(\begin{smallmatrix}a&b&c\\0&0&0\\0&0&0\end{smallmatrix}\right)$ and $\lim_{n\to\infty}\frac{P_nv}{\Vert P_nv\Vert}=\left(\begin{smallmatrix}1\\0\\0\end{smallmatrix}\right)$.

$2^{\rm nd}$ case: $s$ is infinite and $u$ and $\tau$ are finite, in this case the limit-points of $\frac{P_n}{\Vert P_n\Vert}$ have the form $\left(\begin{smallmatrix}a&b&c\\0&0&e\\0&0&f\end{smallmatrix}\right)$. There exists $\lim_{n\to\infty}\frac{P_nv}{\Vert P_nv\Vert}$ independent of $v$ iff $a=b=0$ for all the limit-points, or $c=d=0$ for all the limit-points. In particular, if $\forall n\ b_ne_n\ne0$, then $\lim_{n\to\infty}\frac{P_n}{\Vert P_n\Vert}$ has the form $\left(\begin{smallmatrix}0&0&c\\0&0&e\\0&0&f\end{smallmatrix}\right)$ and $\lim_{n\to\infty}\frac{P_nv}{\Vert P_nv\Vert}=\left(\begin{smallmatrix}c\\e\\f\end{smallmatrix}\right)$.

$3^{\rm rd}$ case: $s$ is finite and $u$ infinite, in this case the limit-points of $\frac{P_n}{\Vert P_n\Vert}$ have the form $\left(\begin{smallmatrix}a&b&c\\0&d&e\\0&0&0\end{smallmatrix}\right)$. There exists $\lim_{n\to\infty}\frac{P_nv}{\Vert P_nv\Vert}$ independent of $v$ iff $a=be-cd=0$ for all the limit-points. In particular, if $\forall n\ b_ne_n\ne0$, then $\lim_{n\to\infty}\frac{P_n}{\Vert P_n\Vert}$ has the form $\left(\begin{smallmatrix}0&b&c\\0&d&e\\0&0&0\end{smallmatrix}\right)$ with $be-cd=0$, and there exists $\lim_{n\to\infty}\frac{P_nv}{\Vert P_nv\Vert}=\left(\begin{smallmatrix}x\\y\\0\end{smallmatrix}\right)$ independent of $v$.

$4^{\rm th}$ case: $s$ and $u$ are finite, in this case the limit-points of $\frac{P_n}{\Vert P_n\Vert}$ have the form $\left(\begin{smallmatrix}a&b&c\\0&d&e\\0&0&f\end{smallmatrix}\right)$. There exists $\lim_{n\to\infty}\frac{P_nv}{\Vert P_nv\Vert}$ independent of $v$ iff $a=d=0$ for all the limit-points, or $d=f=0$ for all the limit-points. In particular, if $\forall n\ b_ne_n\ne0$, then $\lim_{n\to\infty}\frac{P_n}{\Vert P_n\Vert}$ has the form $\left(\begin{smallmatrix}0&0&c\\0&0&e\\0&0&f\end{smallmatrix}\right)$, $\lim_{n\to\infty}\frac{P_nv}{\Vert P_nv\Vert}=\left(\begin{smallmatrix}c\\e\\f\end{smallmatrix}\right)$.

\end{document}